\documentclass[12pt]{article}
\usepackage[margin=1in]{geometry}
\usepackage{amsfonts,amsmath,amssymb}
\usepackage[none]{hyphenat}
\usepackage{fancyhdr}
\usepackage{graphicx}
\usepackage{float}
\usepackage[nottoc,notlot,notlof]{tocbibind}
\usepackage[utf8]{inputenc}
\usepackage{amsthm}
\usepackage{chngcntr}
\usepackage{relsize}
\usepackage{hyperref}
\usepackage{xcolor} 
\hypersetup{ 
    colorlinks,
    linkcolor={black!50!black},
    citecolor={cyan!60!black},
    urlcolor={blue!80!black}
}
\usepackage{listings}
\usepackage{xcolor}
\usepackage{mathrsfs}
\usepackage{tikz-cd}
\usepackage{enumerate}
\usepackage{tablists}
\usepackage{setspace}
\usepackage{scalerel}
\usepackage[title]{appendix}
\usetikzlibrary{shapes,arrows,positioning}

\DeclareFontFamily{U}{mathx}{}
\DeclareFontShape{U}{mathx}{m}{n}{<-> mathx10}{}
\DeclareSymbolFont{mathx}{U}{mathx}{m}{n}
\DeclareMathAccent{\widehat}{0}{mathx}{"70}
\DeclareMathAccent{\widecheck}{0}{mathx}{"71}

\counterwithin{equation}{section}
\theoremstyle{definition}
\newtheorem{theorem}{Theorem}[section]
\newtheorem{defi}[theorem]{Definition}
\newtheorem{prop}[theorem]{Proposition}
\newtheorem{question}[theorem]{Question}
\newtheorem{conjecture}[theorem]{Conjecture}
\newtheorem{corollary}[theorem]{Corollary}
\newtheorem{lema}[theorem]{Lemma}
\newtheoremstyle{exe}
{\topsep}
{\topsep}
{}
{}
{}
{:}
{5pt plus 1pt minus 1pt}
{}
\theoremstyle{exe}
\newtheorem{example}[theorem]{\textnormal{\textit{Example}}}
\newlength\myheight  
\AtEndEnvironment{example}{\null\hfill\ensuremath{\scaleto{\diamond}{\myheight}}}

\newtheoremstyle{note}
{\topsep}
{\topsep}
{}
{}
{}
{:}
{5pt plus 1pt minus 1pt}
{}
\theoremstyle{note}
\newtheorem{remark}[theorem]{\textit{Remark}}
\newenvironment{sproof}{%
  \proof}{\endproof}

\newcolumntype{"}{@{\hskip\tabcolsep\vrule width 2pt\hskip\tabcolsep}}

\newcommand{\Span}{\text{span}}
\newcommand{\im}{\text{im}\,}
\newcommand{\hil}{\mathcal{H}}

\newcommand{\A}{\mathcal{A}}

\newcommand{\Li}{\mathcal{L}}

\newcommand{\cls}{\overline{\text{span}}}
\newcommand{\supp}{\text{supp}}
\newcommand{\nat}{\mathbb{N}}

\newcommand{\id}{\text{id}}

\newcommand{\G}{\mathcal{G}}

\newcommand{\tr}{\text{tr}}
\newcommand{\vN}{\text{vN}}

\newcommand{\Init}{\text{Init}}
\newcommand{\Fin}{\text{Fin}}
\newcommand{\gob}{\mathcal{G}^{(0)}}
\newcommand{\gprod}{\mathcal{G}^{(2)}}
\renewcommand{\L}{L^2(\G,\nu^{-1})}

\newcommand{\ghilb}{(\gob,\{\hil_x\},\mu,L)}
\newcommand{\lrghb}{(\gob,\{\ell^2(\G^x)\},\mu,\lambda)}

\newcommand{\cb}{\Lambda_{\text{cb}}}

\newcommand{\mhotimes}{\otimes_{\text{hX}}}

\usepackage{sectsty}
\sectionfont{\fontsize{15}{15}\selectfont}
\subsectionfont{\fontsize{13}{13}\selectfont}

\newcommand{\ignorar}[1]{}

\makeatletter
\newcommand*\bigcdot{\mathpalette\bigcdot@{.5}}
\newcommand*\bigcdot@[2]{\mathbin{\vcenter{\hbox{\scalebox{#2}{$\m@th#1\bullet$}}}}}
\makeatother

\begin{document}

\title{\sc On weakly amenable groupoids\thanks{Work funded by FCT/Portugal through project UIDB/04459/2025.}}
\author{\sc Tom{\'a}s Pacheco\\
~\\
\begin{minipage}{\textwidth}
\setstretch{1}
\it\normalsize
\center Center for Mathematical Analysis, Geometry and Dynamical Systems,\\
Department of Mathematics, Instituto Superior T\'{e}cnico,
University of Lisbon,\\
Av.\ Rovisco Pais 1, 1049-001 Lisboa, Portugal\\
~\\
E-mail: {\sf tomas.pacheco@tecnico.ulisboa.pt}
\end{minipage}
}
\date{~}
\maketitle
\singlespacing
\setstretch{1}

\begin{abstract}
    In this work, we study groupoids and their approximation properties, generalizing both the definitions and some known results for the group case. More precisely, we introduce weak amenability for groupoids using the definition of the Fourier algebra given by Renault in \cite{renault1997fourier}. We prove that weakly amenable groupoids are inner exact. We also generalize its algebraic counterpart, the CBAP. To do this we introduce the notion of a \textit{quasi Cartan pair} $(B,A)$ and see that $(C_r^*(\G),C_0(\gob))$ can be viewed as such.  We then define what it means for a pair $(B,A)$ to have the CBAP. We introduce the Cowling-Haagerup constants associated to these approximation properites and prove that $\cb(C_r^*(\G),C_0(\gob)) \leq \cb(\G)$. We then study some classes of groupoids where we could achieve equality, that is, $\cb(\G) = \cb(C_r^*(\G),C_0(\G))$. They are discrete groupoids and groupoids arising from partial actions of a discrete group $\Gamma$ on a locally compact Hausdorff space $X$. 
\end{abstract}
\newpage

\tableofcontents
\thispagestyle{empty}
\clearpage

\section{Introduction}

In this work, we study groupoids and their $C^*$-algebras, more specifically, we will study their approximation properties, whose relevance is best explained in the beginning of the book of Brown and Ozawa \cite{brown2008textrm}. "Approximation is ubiquitous in mathematics; when the object in question cannot be studied directly, we approximate by tractable relatives and pass to a limit. In our context this is particularly important because $C^*$-algebras are (almost always) infinite dimensional and we can say precious little without the help of approximation theory. (...) We intend to celebrate it. This subject is a functional analyst's delight, a beautiful mixture of hard and soft analysis, pure joy for the technical inclined.". \\
\\
Let $\Gamma$ be a discrete group and $\lambda: \Gamma \rightarrow B(\ell^2(\Gamma))$ the left regular representation which is given by $\lambda_s\delta_t = \delta_{st}$ for $s,t \in \Gamma$. The \textit{reduced group $C^*$-algebra} is the $C^*$-algebra generated by the image of $\lambda$ in $B(\ell^2(\Gamma))$ and is denoted by $C_r^*(\Gamma)$, moreover, the group von Neumann algebra is the double commutant $\vN(\Gamma) := C_r^*(\Gamma)'' \subseteq B(\ell^2(\Gamma))$. Probably the most important approximation property of them all is amenability of a group $\Gamma$. It was introduced by von Neumann in 1927 to study the Banach-Tarski paradox \cite{anantharaman2001amenable} and in the 1940's it has been essential to the theory of groups and their $C^*$-algebras. Since then, it has penetrated into von Neumann algebras, operator spaces and even differential geometry \cite{runde2004lectures}. However, it was in 1973 in the seminal paper of Lance \cite{lance1973nuclear} where amenability was stated in the form most relevant to this work with the goal of proving that a discrete group $\Gamma$ is amenable if and only if its reduced $C^*$-algebra $C_r^*(\Gamma)$ is nuclear. More concretely, he proved the following theorem.
\begin{theorem}[{\cite[Theorem 3.4]{lance1973nuclear}}]
A $C^*$-algebra $A$ is nuclear if and only if there exists a net of finite rank completely positive operators of norm one $T_i : A \rightarrow A$ such that $T_i \rightarrow \id_A$ in the strong operator topology.   
\end{theorem}
The latter condition of the previous theorem is called the \textit{completely positive approximation property} and is abbreviated as CPAP. Moreover, Lance proved that $\Gamma$ is amenable if and only if $C_r^*(\Gamma)$ has the CPAP \cite[Theorem 4.2]{lance1973nuclear}, obtaining the fact that $\Gamma$ is amenable if and only if $C_r^*(\Gamma)$ is nuclear as a corollary. In this work however, we focus on a weaker notion, introduced by Cowling and Haagerup in \cite{cowling1989completely}, which understandably so is called weak amenability. Moreover, to quantify how strongly this property is satisfied, we associate to the group $\Gamma$ a constant denoted by $\cb(\Gamma)$, called the Cowling-Haagerup constant. Weak amenability is described using the Fourier-Stieltjes $B(\Gamma)$ and the Fourier algebra $A(\Gamma)$ first introduced for abelian groups as $B(\Gamma) = M(\widehat \Gamma)$ and $ A(\Gamma) = L^1(\widehat \Gamma)$, where $\widehat \Gamma$ is the Pontryagin dual and later generalized to all locally compact groups by Eymard in \cite{eymard1964algebre}. $B(\Gamma)$ is defined as the space of coefficients of unitary representations of $\Gamma$, that is, functions of the form
$$
\varphi(t) = (\xi,\eta)_\pi(t) :=  \langle \pi(t)\xi,\eta \rangle_H,
$$
where $\pi: \Gamma \rightarrow B(H)$ is a unitary representation and $\xi,\eta \in H$. For $\varphi \in B(\Gamma)$ we define the norm
$$
\| \varphi \|_{B(\Gamma)} := \inf \{ \|\xi\|_H\|\eta\|_H: \varphi = (\xi,\eta)_\pi \text{ and } \pi: \Gamma \rightarrow B(H) \text{ is a unitary representation} \}.
$$
The Fourier algebra is then the subalgebra $A(\Gamma) \subseteq B(\Gamma)$ of coefficients associated to the left regular representation $\lambda: \Gamma \rightarrow B(\ell^2(\Gamma))$. With the Fourier algebra, we can find yet another equivalent condition to amenability, due to Leptin in \cite{leptin1968algebre}\footnote{In the book of Brown and Ozawa \cite{brown2008textrm} in page 48 it is written that "Amenable groups admit approximately $10^{10^{10}}$ different characterizations".}.
\begin{prop}\label{introleptin}
    A discrete group $\Gamma$ is amenable if and only if there exists a net of finitely supported functions $\varphi_i \in A(\Gamma)$ such that $\varphi_i \xrightarrow{} 1$ pointwise and $\sup_i \|\varphi\|_{B(\Gamma)} < \infty$.
\end{prop}
An important fact about weak amenability, and more concretely, $\cb(\Gamma)$ is that it is an invariant for the group $C^*$-algebra and von Neumann algebra, moreover, in the abstract harmonic analysis setting it provides a uniform bound for approximate units on the Fourier algebra $A(\Gamma)$. The norm where this bound is achieved is not the one coming from $A(\Gamma)$, otherwise $\Gamma$ would be amenable by Proposition \ref{introleptin}, but a different one which we present now since its definition revolves around a phenomenon that is central to this work. If $\varphi \in \ell^\infty(\Gamma)$, it may happen that the ultraweakly continuous extension to $\vN(\Gamma)$ of the map $m_\varphi : C_c(\Gamma) \ni f \rightarrow \varphi f \in C_c(\Gamma)$ is well defined and completely bounded and the space of such functions is denoted by $M_0A(\Gamma)$. It is endowed with a norm defined as $\|\varphi\|_{M_0A(\Gamma)} := \|m_\varphi\|_{\text{cb}}$ and satisfies $\|\varphi\|_{M_0A(\Gamma)} \leq \|\varphi\|_{B(\Gamma)}$. 
\begin{defi}
    A discrete group $\Gamma$ is said to be \textit{weakly amenable} if there exists a net of finitely supported functions $(\varphi_i) \in  C_c(\Gamma)$ and a $C > 0$ such that $\varphi_i \xrightarrow{} 1$ pointwise and $ \sup_i \|\varphi_i\|_{M_0A(\Gamma)} = C < \infty $.
    The infimum of such $C$'s is called the \textit{Cowling-Haagerup constant} and is denoted by $\cb(\Gamma)$. If $\Gamma$ is not weakly amenable, we set $\cb(\Gamma) = + \infty$.
\end{defi}
So far, for a group $\Gamma$ we have introduced its amenability in the $C^*$-algebraic realm, using a net of operators $M_i: C_r^*(\Gamma) \rightarrow C_r^*(\Gamma)$, on the other hand, to describe weak amenability we used functions on $\Gamma$. This interplay is the key phenomenon of this work, indeed there are equivalent definitions for both of these properties that correspond to the ones we already established for the opposed realm. For amenability in terms of functions we have:
\begin{theorem}[{\cite[Theorem 2.6.8]{brown2008textrm}}]
A discrete group $\Gamma$ is amenable if and only if there exists a net of finitely supported positive definite functions $\varphi_i \in \ell^\infty(\Gamma)$ such that $\varphi_i \rightarrow 1$.
\end{theorem}
A positive definite function $\varphi \in \ell^\infty(\Gamma)$ induces by pointwise multiplication a completely positive operator $m_\varphi : C_r^*(\Gamma) \rightarrow C_r^*(\Gamma)$ given by $f \mapsto \varphi f$ and it can be seen that the correspondence $\varphi \mapsto m_\varphi$ proves that if $\Gamma$ is amenable then $C_r^*(\Gamma)$ has the CPAP, this is similar to what Lance did in \cite{lance1973nuclear}. On the other hand, for weak amenability in the $C^*$-world we have:
\begin{theorem}[{\cite[Theorem 2.6]{haagerup2016group}}]\label{introweakamenabilitytheorem}
    A discrete group $\Gamma$ is weakly amenable if and only if there exists a a net of finite rank completely bounded operators $M_i : C_r^*(\Gamma) \rightarrow C_r^*(\Gamma)$ and $C > 0$ such that $M_i \rightarrow \id_{C_r^*(\Gamma)}$ in the strong operator topology and $\sup_i \|M_i\|_{\text{cb}} \leq C$. Moreover, if we denote by $\cb(C_r^*(\Gamma))$ the infimum of such $C$'s we have 
    $$
\cb(\Gamma) = \cb(C_r^*(\Gamma)).
    $$
\end{theorem}
Similarly to the case of amenability, the latter condition of the preceding theorem is called \textit{completely bounded} \textit{approximation property} and is abbreviated by CBAP. 
\begin{sproof}
Every $\varphi \in A(\Gamma)$ induces a completely bounded multiplication operator $m_\varphi: C_r^*(\Gamma) \rightarrow C_r^*(\Gamma)$ such that $\|m_\varphi \|_{\text{cb}} \leq \|\varphi\|_{M_0A(\Gamma)}$, from this it can be obtained that $\cb(C_r^*(\Gamma)) \leq \cb(\Gamma)$. The difficult part is doing the converse. Let $T: C_r^*(\Gamma) \rightarrow C_r^*(\Gamma)$ be a completely bounded operator, the proof can be broken down into the following steps:
\begin{enumerate}
    \item Define a function $\varphi_T \in A(\Gamma)$ associated with $T$ by 
    $$
\varphi_T(t) = \tr(\lambda_t^* T\lambda_t),
    $$
    where $\tr: C_r^*(\Gamma) \rightarrow \mathbb C$ is the canonical faithful trace $\tr(f) = \langle \lambda(f)\delta_e,\delta_e \rangle$.
    \item Obtain a new operator $S:C_r^*(\Gamma) \rightarrow C_r^*(\Gamma)$  such that $S = m_{\varphi_T}$ by applying completely bounded transformations to $T$, more specifically:
    \begin{enumerate}[(i)]
        \item Use Fell absorption to obtain a coaction $\pi:C_r^*(\Gamma) \rightarrow C_r^*(\Gamma) \otimes C_r^*(\Gamma)$, given by $\pi(\lambda_t) = \lambda_t \otimes \lambda_t$.
        \item There exists a conditional expectation $\mathcal{E}: C_r^*(\Gamma) \otimes C_r^*(\Gamma) \rightarrow \im(\pi)$ such that for $s,t \in \Gamma$, $\mathcal{E}(\lambda_s \otimes \lambda_t) = \delta_{s,t} \, \lambda_s$.
        \item The desired operator is $S = \pi^{-1} \circ \mathcal{E} \circ (T \otimes 1) \circ \pi$, meaning that $S = m_{\varphi_T}$ and 
        $$
\|m_{\varphi_T}\|_{\text{cb}} = \|S\|_{\text{cb}} \leq \|T\|_{\text{cb}}.
        $$
        The correspondence $T \mapsto \varphi_T$ will then give $\cb(\Gamma) \leq \cb(C_r^*(\Gamma))$.
    \end{enumerate}
\end{enumerate}
\end{sproof}
The main goal of this work is to generalize the previous theorem in order to apply it to groupoids, although we were not able to achieve this in its strongest form. Groupoids first appeared in operator algebras to provide a link between ergodic theory and von Neumann algebras in the paper of Hahn \cite{hahn1978regular} and later, the $C^*$-algebra setting was developed by Renault in 1978 for his thesis \cite{renault2006groupoid}. Due to its relevance for groups, it was not much later that amenability was introduced for groupoids, its first instances were done in the end of the $70$'s and are due to Zimmer in \cite{zimmer1977hyperfinite, zimmer1977neumann, zimmer1978amenable} where he defined amenability for discrete group actions and countable equivalence relations. In the following years it was generalized to all sorts of contexts playing a big role in many advances in certain topics in operator algebras: exact $C^*$-algebras, Elliott's program, Baum-Connes conjecture and Novikov conjecture, to name a few. More information on this can be found in the preface or introduction of the monograph \cite{anantharaman2001amenable}. As we just mentioned, amenability is widely present in groupoid literature, however there is no account on weak amenability to be found, the closest being the work of Mckee in \cite{mckee2016weak} where weak amenability was generalized to $C^*$-dynamical systems. The main goal of this work is to fill in this gap. The first step to define weak amenability is to introduce the Fourier algebra of a groupoid, luckily, this has already been done by Renault in \cite{renault1997fourier}. To do this, it is essential to study the representation theory of groupoids, which is considerably more involved than the group case. The object of the representations is a measured Hilbert bundle over the unit space $\gob$ of the groupoid $\G$, it is a triple $\hil := (\gob, \{ \hil_x\}, \mu)$ where $\mu$ is a Borel probability measure on $X$ and for each $x \in X$, $\hil_x$ is the fiber at $x$ which is a Hilbert space. Representations are thus given by a family of unitary operators $L(\gamma): B(\hil_{d(\gamma)}) \rightarrow B(\hil_{r(\gamma)})$ indexed by elements $\gamma \in \G$ and that are compatible with the groupoid operations. Coefficients of such representations are defined using two bounded sections $\xi,\eta \in L^\infty(X,\hil)$ and letting 
$$
(\xi,\eta)(\gamma) = \langle \xi \circ r(\gamma), L(\gamma)(\eta \circ d(\gamma)) \rangle_{\hil_{r(\gamma)}}, \quad \text{ for all } \gamma \in \G.
$$
The Fourier-Stieltjes algebra will be the space of coefficients of representations of $\G$ equipped with the norm
$$
\| \varphi \|_{B(\G)} := \inf \{ \|\xi\|_\infty \|\eta \|_\infty: \varphi = (\xi,\eta)_L \text{ for a representation } \ghilb \}.
$$
\indent As with the group case, there is a distinguished representation, the left regular representation. For groupoids, it is given by the bundle $\Li^2_\mu(\G) := (\gob, \ell^2(\G^x),\mu)$ where $\mu$ is a quasi-invariant measure on $\gob$ and the unitary operators $\lambda_\gamma \delta_\beta = \delta_{\gamma\beta}$ for $\gamma \in \G$ and $\beta \in \G^{d(\gamma)}$. The Fourier algebra $A(\G)$ of a measured groupoid $\G$ is then defined as the closure in $B(\G)$ of the coefficients associated to the left regular representation (see \cite{renault1997fourier}). When $\G = \Gamma$ is a discrete group, there was no need to take the closure since $A(\Gamma)$ is already closed in $B(\Gamma)$, however for groupoids this is no longer the case. Moreover, some authors decided to not take closures, for example in \cite{paterson2004fourier}. With this, the definition of weak amenability generalizes to measured groupoids in a similar fashion as it was done for amenability in \cite{anantharaman2001amenable}: 
\begin{defi}
    A  measured groupoid $(\G,\mu)$ is said to be \textit{weakly amenable} if there exists a net $(\varphi_i) \in A(\G)$ and a $C>0$ such that $\sup_i \|\varphi_i \|_{M_0A(\G)} = C < \infty$ and $\varphi_i \xrightarrow{} 1$ in the weak-$*$ topology of $L^\infty(\G)$. The \textit{Cowling-Haagerup constant} is the infimum of all $C > 0$ for which there is a net as above and is denoted by $\Lambda_{\text{cb}}(\G,\mu)$.
\end{defi}
We introduce also topological definitions of weak amenability.
\begin{defi}
    Let $\G$ be an étale groupoid, then $\G$ is said to be \textit{measurewise weakly amenable} if there exists a $C > 0$ such that for all quasi-invariant probability measures $\mu$ on $\gob$, $(\G,\mu)$ is weakly amenable with $\cb(\G,\mu) \leq C$. The \textit{measurewise Cowling-Haagerup constant} is the infimum of all $C > 0$ for which above holds and is denoted by $\Lambda_{\text{mcb}}(\G)$. 
\end{defi}
For topological amenability, we will want to make use of continuous functions and since $\G$ is étale, we have that $C_c(\G) \subseteq A(\G)$, where we endow $\G$ with a quasi-invariant measure with full support.
\begin{defi}
    Let $\G$ be an étale groupoid, then $\G$ is said to be \textit{(topologically) weakly amenable} if there exists  a net $(\varphi_i) \in C_c(\G)$, a $C > 0$ and quasi-invariant probability measure with full support $\mu$, such that $\sup_i \|\varphi_i \|_{M_0A(\G)} = C < \infty$ and $\varphi_i \xrightarrow{} 1$ uniformly on compact subsets. The \textit{Cowling-Haagerup constant} is the infimum of all $C > 0$ for which above holds and is denoted by $\Lambda_{\text{cb}}(\G)$.
\end{defi}
We remark thatt the value of $\Lambda_{\text{cb}}(\G)$ does not depend on the chosen quasi-invariant measure with full support. From the definitions, it can be immediately seen that for measured groupoids, amenability implies weak amenability and the same holds for étale groupoids, both in the measurewise and topological case. In this work, we will also prove that if $\G$ is topologically weakly amenable, then it is inner exact. Thus, the work of Willet in \cite{willett2015non} using the HLS groupoids of Higson, Lafforgue and Skandalis allows us to say that there exist groupoids that have the weak containment property but are not weakly amenable. On the other hand, all non-amenable weakly amenable groups do not have weak containment.\\
\\
We can associate certain operator algebras to a measured groupoid $(\G,\mu)$. Using the left regular representation we can construct the reduced $C^*$-algebra $C_r^*(\G,\mu) \subseteq B(\L)$ and the von Neumann algebra $\vN(\G,\mu) :=  C_r^*(\G,\mu)'' \subseteq B(\L)$. There also exists a \textit{full} $C^*$-algebra, denoted by $C^*(\G,\mu)$. We remark that the definition of $C_r^*(\G,\mu)$ depends only on the support of $\mu$ \cite[Section 1.3]{anantharaman2016some} , however $\vN(\G,\mu)$ and $C^*(\G,\mu)$ depend on the measure class of $\mu$ rather than just the support \cite[p.5]{anantharaman2013haagerup}. These algebras can be related to the Fourier-Stieltjes and Fourier algebras via a duality theory similar to what Eymard did in \cite{eymard1964algebre} for (discrete) groups. Indeed, if $\G = \Gamma$ is a discrete group, then $B(\Gamma) \cong C^*(\Gamma)^*$ and $A(\Gamma) \cong \vN(\Gamma)_*$. In his paper \cite{renault1997fourier}, Renault obtained a generalization of this fact by defining the spaces
$$
\mathcal{X}(\G) := L^2(X)^* \mhotimes C^*(\G,\mu) \mhotimes L^2(X)
$$
and 
$$
\mathcal{Y}(\G) := L^2(X)^* \mhotimes A(\G) \mhotimes L^2(X)
$$
where $\mhotimes$ is the module Haagerup tensor product over $A := L^\infty(\gob)$. Renault proved that $B(\G) \cong \mathcal{X}(\G)^*$ and $\vN(\G)_* \cong \mathcal{Y}(\G)$. Let $\varphi \in L^\infty(\G)$ and $m_\varphi: \vN(\G,\mu) \rightarrow \vN(\G,\mu)$ be the induced multiplication map $f \mapsto \varphi f$. \\
\\
We now introduce the corresponding notion on the $C^*$-level to weak amenability, generalizing the CBAP to the groupoid context. The fact we are using a $C^*$-algebraic theory suggest that the more tractable variant of weak amenability will be the topological one and indeed, it will be the focus of this work. In the case where $m_\varphi$ is well defined and bounded, $\varphi$ is said to be a \textit{bounded multiplier} and we write $\varphi \in MA(\G)$. If in addition, $m_\varphi$ is completely bounded, then $\varphi$ is a \textit{completely bounded multiplier} and we write $\varphi \in M_0A(\G)$. Looking at the proof of Theorem \ref{introweakamenabilitytheorem}, to generalize the completely bounded approximation property (CBAP) to groupoids, we must ask ourselves what properties does $m_\varphi$ have when $\varphi \in C_c(\G)$. Of course, we should expect that $m_\varphi$ is completely bounded and this is true since $C_c(\G) \subseteq M_0A(\G)$ \cite[Proposition 3.3]{renault1997fourier}. However, we cannot guarantee that $m_\varphi$ will have finite rank, for example, if we let $X$ be a compact Hausdorff space endowed with a Borel probability measure $\mu$ with full support then $\G = \gob = X$ is a measured groupoid and $C_r^*(X) = C(X)$. Then, the multiplier induced by the constant $1$ function $m_1 = \id_{C(X)}$ does not have finite rank. To solve this issue, we borrow inspiration from \cite{anantharaman2013haagerup}, where the Haagerup property was generalized to groupoids and notice that if $\varphi \in C_c(\G)$, then $m_\varphi$ is not only linear but also $C_0(\gob)$-linear with action given by
$$
f * a(\gamma) = f(\gamma) \, a(s(\gamma)), \text{ for all } f \in C_c(\G) \text{ and } a \in C_0(\gob).
$$
So instead of requiring that $m_\varphi$ has finite rank, we demand finite $C_0(\gob)$-rank. To explain what this means more precisely we make use of the framework of theory of Hilbert $C^*$-modules. If $\G$ is étale and Hausdorff, then there exists a conditional expectation $E:C_r^*(\G) \rightarrow C_0(\gob)$ extending the restriction map $C_c(\G) \rightarrow C_0(\gob)$ \cite[p.204]{brown2008textrm} which allows us to define a $C_0(\gob)$-valued inner product by 
$$
\langle f,g \rangle_{C_0(\gob)} := E(f^* * g).
$$
Together with the above action, this turns $C_r^*(\G)$ into a pre-Hilbert $C^*$-module \cite[p.203]{brown2008textrm}. The notion of rank one operators for Hilbert $C^*$-modules is for example defined in \cite[II.7.2.5]{blackadar2006operator}, which for our case is as follows: an operator $T:C_r^*(\G) \rightarrow C_r^*(\G)$ is said to have \textit{$C_0(\gob)$-rank one} if there are $g,h \in C_r^*(\G)$ such that for all $f \in C_r^*(\G)$
$$
Tf = h * \langle g,f \rangle_{C_0(\gob)} = h * E(g^* * f).
$$
With this in mind, a finite $C_0(\gob)$-rank operator is one that is a finite linear combination of $C_0(\gob)$-rank one operators and indeed, in Proposition \ref{finiterankmultipliers}, we prove that if $\varphi \in C_c(\G)$, then $m_\varphi$ has finite $C_0(\gob)$-rank. We can almost fit our setting in the context of Cartan pairs as in \cite{renault2024cartan}.
\begin{defi}
    Let $A \subseteq B$ be an abelian sub $C^*$-algebra of a $C^*$-algebra $B$. We will say that $A$ is a \textit{Cartan subalgebra} if:
    \begin{enumerate}[(i)]
    \itemsep-5pt
        \item $A$ contains an approximate unit of $B$;
        \item $A$ is regular;
        \item there exists a conditional expectation $E:B \rightarrow A$. 
    \end{enumerate}
    Then $(B,A)$ is called a \textit{quasi Cartan pair}. If moreover 
    \begin{enumerate}[(iv)]
        \item $A$ is maximal abelian.
    \end{enumerate}
    Then $(B,A)$ is called a \textit{Cartan pair}.
\end{defi}
If $\G$ is an étale groupoid, then $(C_r^*(\G),C_0(\gob))$ is a quasi Cartan pair, additionally, if $\G$ is essentially principal, then $(C_r^*(\G),C_0(\gob))$ is a Cartan pair. In this work, we will not need condition (iv) so we will focus on quasi Cartan pairs.
\begin{defi}
  A quasi Cartan pair $(B,A)$ is has the \textit{completely bounded approximation property} (abbreviated by CBAP) if there exists a net $T_i: B \rightarrow B$ of $A$-linear operators with finite $A$-rank such that $T_i \xrightarrow{\text{SOT}} \id_B$ and $\sup_i \| T_i \|_{\text{cb}} = C < \infty$. The \textit{Cowling-Haagerup constant} is the least $C$ such that there is a net as above and is denoted by $\Lambda_{\text{cb}}(B,A)$. 
\end{defi}
Note that if $\G = \Gamma$ is a group, then $\gob = \{e\}$ and so $C_0(\gob) = \mathbb C$ and the CBAP of $(C_r^*(\G),C_0(\gob))$ reduces to the normal CBAP.  We would like to prove that for a groupoid $\G$, we have 
$$
\cb(\G) = \cb(C_r^*(\G),C_0(\gob)).
$$
however, we were not able to show this in general and only for some subclasses of groupoids (see Chapter 5). If $\varphi \in C_c(\G)$, the correspondence $\varphi \mapsto m_\varphi$ yields 
$$
 \cb(C_r^*(\G),C_0(\gob)) \leq \cb(\G), 
$$
For the converse, let $T: C_r^*(\G) \rightarrow C_r^*(\G)$ and define
$$
\varphi_T(\gamma) = E( \chi_s^* * T\chi_s) \circ d(\gamma)
$$
where $s$ is any compact open bisection such that $\gamma \in s$. At this point, we assume that $\G$ has a basis of such bisections, that is, $\G$ is \textit{ample}, so it can be seen that $\varphi_T$ is well defined, i.e. does not depend on the choice of $s$ and it is worth to take a moment and compare this with the definition of $\varphi_T$ in the proof of Theorem \ref{introweakamenabilitytheorem} for the group case. The main obstruction was given by the fact that we could not prove that $\|\varphi\|_{M_0A(\G)} \leq \|T\|_{\text{cb}}$, which is caused by the fact that a Fell absorption principle that yields a "diagonal" coaction $C_r^*(\G) \rightarrow C_r^*(\G) \otimes \A$ for a suitable $\A$ as in Theorem \ref{introweakamenabilitytheorem} is not available. However, for some relevant subclasses of groupoids it exists, such as discrete groupoids and groupoids arising from partial actions of a discrete group $\Gamma$ on a locally compact Hausdorff space $X$. In these cases, we employ the appropriate version of Fell's absorption principle and a method similar to the proof of \ref{introweakamenabilitytheorem} to prove if $\G$ is a groupoid of the above classes then 
$$
\Lambda_{\text{cb}}(C_r^*(\G),C_0(\gob)) = \Lambda_{\text{cb}}(\G).
$$
This gives good evidence of the fact that the above equality holds in all cases, which is our desire. \\
\\
This work is structured as follows. In Section 2, we setup the basic theory necessary to delve deeper into approximation properties. We define a groupoid, its topological and measure-theoretic properties, the associated operator algebras $C_r^*(\G),\vN(\G)$ and study its representation theory. In Section 3, we define the Fourier-Stieltjes and Fourier algebras of $\G$ and we follow it by study its multipliers and duality theory in Section 4. We then introduce the central concept of this work, being weak amenability. In Section 5, we define the three different variants: 
\begin{enumerate}[(i)]
    \itemsep-5pt
    \item weak amenability for measured groupoid $(\G,\mu)$;
    \item measurewise weak amenability for étale groupoids $\G$;
    \item topological weak amenability for étale groupoids $\G$.
\end{enumerate}
and study their relationships between each other and with other approximation properties. In particular, we will show:
\begin{enumerate}[(i)]
    \itemsep-5pt
    \item If $(\G,\mu)$ is amenable then it is weakly amenable and $\cb(\G,\mu) = 1$;
    \item If $\G$ is measurewise amenable, then it is measurewise weakly amenable and $\Lambda_{\text{mcb}}(\G) = 1$;
    \item If $\G$ is topologically amenable, then it is topologically weakly amenable and $\cb(\G) = 1$;
    \item If $\G$ is topologically weakly amenable, then it is inner exact;
    \item If $\G$ is topologically weakly amenable, then it is measurewise weakly amenable and $\Lambda_{\text{mcb}}(\G) \leq \cb(\G)$;
\end{enumerate}
In Section 6, we introduce (quasi) Cartan pairs and the completely bounded approximation property associated to them. We finish that section by proving that if $\G$ is an étale groupoid, then
$$
\cb(C_r^*(\G),C_0(\gob)) \leq \cb(\G).
$$
We finish this work with some classes of groupoids of groupoids for which we were able to obtain an equality 
$$
\cb(C_r^*(\G),C_0(\gob)) = \cb(\G).
$$
As we mentioned before, these are groupoids arising from a partial action of a discrete group $\Gamma$ on a locally compact Hausdorff space $X$, studied in Section 7 and discrete groupoids in Section 8. \\
\\
\textbf{Conventions:} We will assume that all topological groupoids $\G$ are locally compact second countable spaces with locally compact Hausdorff unit space $X := \gob$. This implies that the fibers $\G^x$ and $\G_x$ are countable.\\
\\
\textbf{Acknowledgments:} The author wishes to thank his thesis supervisor P. Resende, for letting this work be possible and Professor Alcides Buss for his immense help. The author also wishes to thank the Center for Mathematical Analysis, Geometry and Dynamical Systems (CAMGSD) and Fundação para a Ciência e Tecnologia (FCT) for the financial support. 
 
\section{Measured Groupoids and their Operator Algebras}
Let $\G$ be a groupoid, that is, a small category in which every arrow is invertible. The set of objects of $\G$ is denoted by $X := \gob$ and we will refer to them as \textit{units} and $X$ as the \textit{unit space}. The structure maps $d,r: \G \rightarrow X$ are called the source and range map, respectively. We endow $\G$ with a topology such that $d,r$ are continuous, $\G$ is a locally compact space and $X$ is Hausdorff (in addition to being locally compact). When this is the case, $\G$ is said to be a \textit{locally compact groupoid}, furthermore, for each $x \in X$, we will want to introduce a measure on each of the fibers $\G_x := d^{-1}(x),\G^x := r^{-1}(x)$ that will be in some way invariant for the groupoid operation; in all its generality this is called a Haar system \cite[Definition 2.2]{renault2006groupoid} which mimics the Haar measure construction for a group. However, to avoid complications, we will restrict to the case where $\G$ is \textit{étale}, meaning the above mentioned fibers are discrete subspaces of $\G$ and we endow the counting measure on them. By a \textit{bisection} we mean a subset $U \subseteq \G$ such that the structure maps $d|_U,r|_U$ are homeomorphisms and we denote the set of all open bisections by $\G^{op}$, moreover, we denote the set of all compact open bisections by $\G^a$. It is known that $\G$ is étale if and only if $\G^{op}$ is a basis for the topology of $\G$ \cite[p.44]{paterson2012groupoids} and we will say that $\G$ is \textit{ample} if $\G^a$ is also a basis for the topology of $\G$. There is a last piece of structure we introduce on $\G$, a Borel probability measure $\mu$ on $X$ that via the counting measure on each fiber induces
$$
\nu(A) = \int_X |A \cap \G^x| \, d\mu(x),
$$
which is a measure on the whole of $\G$, moreover, if we define $\nu^{-1}$ in the same way as above replacing $\G^x$ with $\G_x$ we get a new measure denoted by $\nu^{-1}$ and $\mu$ is said to be \textit{quasi-invariant} if $\nu$ and $\nu^{-1}$ have the same null sets. When $\G$ is endowed with such a $\mu$, $\G$ is said to be a \textit{measured groupoid}\footnote{The usual definition of measured groupoids does not presuppose a topology, but in our case, we will assume that a measured groupoid is an \textbf{étale} groupoid endowed with a quasi-invariant measure $\mu$ and Haar system given by the counting measure.}. Quasi-invariant measures always exist, see \cite[Remark 3.18]{muhly1997coordinates}.\\
\\
\textbf{$C^*$-algebras of Hausdorff étale groupoids.}
Let $\G$ be a Hausdorff étale groupoid with unit space $X:= \gob$ and $C_c(\G)$ be the set of continuous functions $f:\G \xrightarrow{} \mathbb C$ with compact support. We endow this space with a convolution product in the following manner: for $\gamma \in \G$ and $f,g \in C_c(\G)$,
$$
f*g(\gamma) = \sum\limits_{\alpha\beta = \gamma}f(\alpha)g(\beta) = \sum\limits_{\beta \in \G_{d(\gamma)}}f(\gamma\beta^{-1})g(\beta),
$$
and an involution $f^*(\gamma) = \overline{f(\gamma^{-1})}$. Define the \textit{$I$-norm} for $f \in C_c(\G)$ by 
$$
\|f\|_I := \max \left\{ \sup\limits_{x \in X} \sum\limits_{\gamma \in \G^x} |f(\gamma)|\, , \, \, \sup\limits_{x \in X} \sum\limits_{\gamma \in \G_x} |f(\gamma)| \right\}.
$$
A \textit{representation} of $\G$ is an $I$-norm continuous homomorphism $C_c(\G) \rightarrow B(H)$. The universal norm on $C_c(\G)$ is thus 
$$
\|f\|_* = \sup\{ \|\pi(f)\|: \pi \text{ is a representation of } C_c(\G) \},
$$
and the \textit{full $C^*$-algebra} of $\G$, denoted by $C^*(\G)$ is the completion of $C_c(\G)$ under $\|\cdot\|_*$. To define the reduced version, we endow $C_c(\G)$ with a right $C_0(X)$-module structure by letting
$$
f  \cdot a(\gamma) = f * a(\gamma) = f(\gamma)a(d(\gamma)), \quad \text{ for } f \in C_c(\G) \text{ and } a \in C_0(X),
$$
and to finish a $C_0(X)$-valued inner product defined for $f,g \in C_c(\G)$ and $x \in X$ as
\begin{equation}\label{moduleinnerproduct}
   \langle f,g \rangle(x) = \sum\limits_{\gamma \in \G_x}\overline{f(\gamma)}g(\gamma). 
\end{equation}
Denote by $L^2(\G)$ the Hilbert $C_0(X)$-module obtained through the completion of $C_c(\G)$ under the norm induced by $\langle \cdot,\cdot \rangle$ and define the \textit{left regular representation} $\lambda: C_c(\G) \xrightarrow{} B(L^2(\G))$  
$$
\lambda(f)g = f*g, \quad \text{ for } f,g \in C_c(\G).
$$
The \textit{reduced groupoid} $C^*$\textit{-algebra} is therefore the completion of $C_c(\G)$ with respect to the norm $\|f\|_r := \|\lambda(f)\|$ and is denoted by $C_r^*(\G)$.\\
\\ 
\textbf{Representation theory of measured groupoids.}
A \textit{representation of $\G$ (or a $\G$-Hilbert bundle)} is a tuple $(\hil,L)$, composed of a measurable Hilbert bundle $\hil$ over $(X,\mu)$ and a family of unitary operators $L=\{L(\gamma):\hil_{d(\gamma)} \rightarrow \hil_{r(\gamma)}\}_{\gamma \in \G}$ such that: 
\begin{enumerate}[(i)]
\itemsep0pt
    \item $L(x)$ is the identity map on $\hil_x$, for all $x \in X$;
    \item $L(\beta)L(\gamma) = L(\beta\gamma)$ for $\nu^2$-almost all $(\beta,\gamma) \in \gprod$;
    \item $L(\gamma)^{-1} = L(\gamma^{-1})$ for $\nu$-almost all $\gamma \in \G$;
    \item for any $\xi,\eta \in L^2(X,\hil)$ the function
    $
\gamma \mapsto \langle \xi \circ r(\gamma), L(\gamma)(\eta \circ r(\gamma)) \rangle
    $
    is $\nu$-measurable on $\G$.
\end{enumerate}
Given a $\G$-Hilbert bundle $(\hil,L)$, we can induce a representation $\pi_L$ of $C_c(\G)$ on the space of sections $L^2(X,\hil)$ by 
\begin{equation}\label{inducedrepoid}
\langle \pi_L(f)\xi,\eta \rangle = \int_\G f(\gamma) \, \langle \xi \circ r(\gamma), L(\gamma)(\eta \circ r(\gamma)) \rangle \,  d\nu_0(\gamma).    
\end{equation}
\begin{theorem}[{\cite[Theorem 3.1.1]{paterson2012groupoids}}]
    Let $\G$ be an étale Hausdorff groupoid, then every $I$-norm continuous representation of $C_c(\G)$ is of the form $\pi_L$ for some $\G$-Hilbert bundle $(\hil,L)$ over $(X,\mu)$ for a quasi-invariant measure $\mu$.
\end{theorem}
\begin{example}\label{leftregulargbundle}
We consider the bundle $\Li^2(\G) := (X,\{\ell^2(\G^x)\},\mu)$ and define the \textit{left regular $\G$-Hilbert bundle} by assigning to $\gamma \in \G$ the operator $(\lambda_\gamma f)(\eta) = f(\gamma^{-1}\eta)$ for each $f \in \ell^2(G^{d(\gamma)})$ and $\eta \in \G^{r(\gamma)}$. By equation \ref{inducedrepoid}, this induces a representation of $C_c(\G)$ on the space of sections $L^2(X,\Li^2(\G))$. We can identify $L^2(X,\Li^2(\G))$ with $L^2(\G,\nu)$ as follows: let $\xi \in L^2(X,\Li^2(\G))$ then we define 
$$
\widehat \xi(\gamma) := \xi(r(\gamma))(\gamma) \in L^2(\G,\nu)
$$ 
and this correspondence is easily seen to be unitary. Moreover, the map $f \mapsto D^{.1/2}f$ is a unitary between $L^2(\G,\nu)$ and $\L$, meaning we can identify $L^2(X,\Li^2(\G))$ with $\L$. If we do this, the induced representation of $C_c(\G)$ on $\L$ becomes $\lambda(f)g = f*g$, for $g \in L^2(\G,\nu^{-1})$ (see \cite[Equation 3.36]{paterson2012groupoids}).    
\end{example}
\textbf{Reductions.} Let $Y \subseteq X := \gob$, the subgroupoid $\G|_Y := r^{-1}(Y) \cap d^{-1}(Y)$ is called the \textit{reduction} of $\G$ by $Y$, it is étale Hausdorff if $\G$ also is. Let now $\mu$ be a quasi-invariant measure on $X$ and $\Omega := \supp(\mu)$, it is a closed $\G$-invariant subset, meaning that for all $\gamma \in \G$, $d(\gamma) \in \Omega$ if and only if $r(\gamma) \in \Omega$. Letting $U = X \backslash \Omega$, we have that the inclusion $C_c(\G|_U) \hookrightarrow C_c(\G)$ induces injective $*$-homomorphisms $C^*(\G|_\Omega) \rightarrow C^*(\G)$ and $C^*_r(\G|_\Omega) \rightarrow C^*_r(\G)$. Moreover, the restriction $C_c(\G) \rightarrow C_c(\Omega)$ induces surjective maps $C^*(\G) \rightarrow C^*(\G|_\Omega)$, $C^*_r(\G) \rightarrow C^*_r(\G|_\Omega)$ and the following sequence is exact
\begin{equation}\label{exactreduction}
0 \longrightarrow C^*(\G|_U) \longrightarrow C^*(\G) \longrightarrow C^*(\G|_\Omega) \longrightarrow 0.    
\end{equation}
The same might not happen for the reduced $C^*$-algebra, which can be verified using the infamous HLS groupoid. See \cite[Section 1.3]{anantharaman2016some} for all of the above.
\\
\\
\textbf{Operator algebras of Hausdorff measured groupoids.} Let $(\G,\mu)$ be a Hausdorff measured (étale) groupoid, to define the full $C^*$-algebra we make use of representations on bundles over $(X,\mu)$. Indeed, let $f \in C_c(\G)$ and define
$$
\|f\|_\mu := \sup\{ \|\pi_L(f)\| : (\hil,L) \text{ is a $\G$-Hilbert bundle over } (\G,\mu) \}.
$$
It is a semi-norm on $\G$ and we let $N_\mu := \{f \in C_c(\G): \|f\|_\mu = 0 \}.$ The \textit{full $C^*$-algebra} of $(\G,\mu)$, denoted by $C^*(\G,\mu)$ is the completion of $C_c(\G)/N_\mu$ under $\| \cdot \|_\mu$. Let now $\lambda: C_c(\G) \rightarrow B(\L)$ be the representation of Example \ref{leftregulargbundle}, the \textit{reduced $C^*$-algebra} of $\G$ is defined to be the closure of $\lambda(C_c(\G))$ in $B(\L)$ and denoted by $C^*_r(\G,\mu)$. Moreover, the groupoid von Neumann algebra is defined as $\vN(\G) := C_r^*(\G)'' \subseteq B(\L)$. $C^*$-algebras of reductions are closely related to the ones we just defined, indeed if we let $\Omega := \supp(\mu)$ then $C^*_r(\G|_\Omega) \cong C_r^*(\G,\mu)$ and there exists a surjective $*$-homomorphism $C^*(\G|_\Omega) \rightarrow C^*(\G,\mu)$, although it might not be a isomorphism, again taking the HLS groupoid as an example (see \cite[Section 1.3]{anantharaman2016some}).
\section{Fourier-Stieltjes and Fourier Algebras}
In this section we define the Fourier and Fourier-Stieltjes algebras of measured groupoids introduced by Renault in \cite{renault1997fourier}. We fix a Hausdorff étale measured groupoid $(\G,\mu)$ throughout with unit space $X := \gob$. Let $\nu,\nu^{-1}$ be the measures on $\G$ induced by $\mu$, since $\mu$ is quasi-invariant, $L^\infty(\G,\nu)$ is isometrically isomorphic to $L^\infty(\G,\nu^{-1})$ and with this in mind, we will write both of these spaces as $L^\infty(\G,\mu)$ or simply $L^\infty(\G)$, when $\mu$ is clear from the context.
\begin{prop}[{\cite[Proposition 1.1]{renault1997fourier}}]\label{posdefoid}
Let $\varphi \in L^\infty(\G)$, the following are equivalent 
\begin{enumerate}[(i)]
\itemsep0pt
    \item For every $N \in \nat$ and $\zeta_1,...,\zeta_n \in \mathbb C$,
    $$
\sum\limits_{i=1}^N \varphi(\gamma_i\gamma_j^{-1})\overline{\zeta_i}\zeta_j \geq 0,
    $$
    for almost all $x \in X$ and all $\gamma_1,...,\gamma_n \in \G_x$;
    \item There exists a $\G$-Hilbert bundle $(\hil.L)$ and a section $\xi\in L^\infty(X,\hil)$ such that for all $\gamma \in \G$,
$$
(\xi,\xi)_L(\gamma) := \langle \xi \circ r (\gamma), L(\gamma) \xi \circ s(\gamma) \rangle = \varphi(\gamma).
$$
\end{enumerate}
When the above conditions hold $\varphi$ is said to be \textit{positive definite} and we write $\varphi \in P(\G)$. Moreover, if $\xi \in L^2(X,\hil)$ is any section such that $\varphi = (\xi,\xi)_L$ we have $\|\varphi\|_\infty = \|\xi\|_\infty^2$.
\end{prop}
\begin{proof}
    The fact that (i) and (ii) are equivalent is proven in \cite[Proposition 1.1]{renault1997fourier}, we will only prove the last part. Suppose $\varphi$ is positive definite and $\xi \in L^2(X,\hil)$ is any section with $\varphi = (\xi,\xi)_L$. By Cauchy-Schwarz it is immediate that $\|\varphi\|_\infty \leq \|\xi\|_\infty^2$, on the other hand, if $x \in X$,
    $$
\varphi(x) = (\xi,\xi)_L(x) = |\xi(x)|^2,
    $$
    doing the supremum over all $x \in X$ we obtain that $\|\varphi\|_\infty \geq \|\xi\|_\infty^2$, finishing the proof.
\end{proof}
\begin{prop}[{\cite[Proposition 1.3]{renault1997fourier}}]
Let $\varphi \in L^\infty(\G)$, the following are equivalent 
\begin{enumerate}[(i)]
\itemsep0pt
    \item $\varphi$ is a linear combination of elements in $P(\G)$.
    \item There exists a $\G$-Hilbert bundle $(\hil,L)$ and sections $\xi,\eta \in L^\infty(X,\hil)$ such that for all $\gamma \in \G$,
$$
(\xi,\eta)_L(\gamma) := \langle \xi \circ r (\gamma), L(\gamma) \eta \circ d(\gamma) \rangle = \varphi(\gamma).
$$
\end{enumerate}
\end{prop}
\begin{defi}
Given a $\G$-Hilbert bundle $(\hil,L)$ and sections $\xi,\eta \in L^\infty(X,\hil)$, the function $(\xi,\eta)_L$ defined above is called a \textit{coefficient}.  Elements $\varphi \in L^\infty(\G)$ that verify the above conditions will also be called coefficients. The set of all coefficients is denoted by $B(\G,\mu)$ or $B(\G)$ when $\mu$ is clear from the context.
\end{defi}
Note that $B(\G)$ is clearly an involutive $\mathbb C$-algebra. The fact that is closed under sums reflects taking direct sums of representations, moreover, it is closed under multiplication which reflects the tensor product of  representations. We can now define a norm by 
$$
\| \varphi \|_{B(\G)} := \inf \{ \|\xi\|_\infty \|\eta \|_\infty: \varphi = (\xi,\eta)_L \text{ for a $\G$-Hilbert bundle } (\hil,L) \text{ over } (X,\mu) \}.
$$
\begin{prop}[{\cite[Proposition 1.4]{renault1997fourier}}]
  $B(\G)$ equipped with the norm  $\| \cdot \|_{B(\G)}$ is a Banach $*$-algebra.
\end{prop}
Let us now consider the left regular $\G$-Hilbert bundle $(\Li^2_\mu(\G),\lambda) = \lrghb$ of Example \ref{leftregulargbundle} and the set 
$$
\A(\G,\mu) := \{ \varphi \in L^\infty(\G): \varphi = (\xi,\eta)_\lambda \}
$$
of its coefficients. We define \textit{Fourier algebra of} $\G$ to be the closure of $\A(\G,\mu)$ in $B(\G,\mu)$, that is, 
$$
A(\G,\mu)  := \overline{ \{ \varphi \in L^\infty(\G): \varphi = (\xi,\eta)_\lambda \} } \subseteq B(\G,\mu).
$$
Of course, we can simply write $A(\G)$ when $\mu$ is clear from the context. Although elements of $A(\G)$ might not be coefficients of the left regular Hilbert bundle, we can construct a new one that will allows us view elements of $A(\G)$ as coefficients. Let $\Li^2_\mu(\G,\ell^2) := (X,\{ \ell^2(\G^x, \ell^2) \}, \mu )$, where 
$$
\ell^2(\G^x,\ell^2) = \{ f: \G^x \rightarrow \ell^2 \text{ measurable}: \|f\|_2 <\infty \}.
$$
For $f \in \ell^2(\G^x,\ell^2)$ the norm $\|f\|_2$ is defined as follows: for every $n$, we have a function $f_n: \G^x \rightarrow \mathbb C$ given by $f_n(\gamma) = f(\gamma)_n$, then 
$$
\|f\|_2 := \left( \sum\limits_{n} \|f_n\|_2^2 \right)^{1/2}.
$$
From the $\G$-Hilbert bundle $(\Li^2_\mu(\G),\lambda)$ we can construct a new representation $\lambda^\infty$ of $\G$ on $\Li^2_\mu(\G,\ell^2)$; for $\gamma \in \G$ and $f = (f_n) \in \ell^2(\G^{s(\gamma)},\ell^2)$ we define
$$
\lambda^\infty_\gamma(f) = (\lambda_\gamma f_n) \in \ell^2(\G^{r(\gamma)},\ell^2).
$$
The $\G$-Hilbert bundle $(\Li^2_\mu(\G,\ell^2),\lambda^\infty )$ is called the \textit{left regular $\G$-Hilbert bundle with infinite multiplicity}.
\begin{prop}[{\cite[Lemma 1.2 and Proposition 1.5]{renault1997fourier}}]
    Every element $\varphi \in A(\G)$ can be written as a coefficient of the left regular $\G$-Hilbert bundle with infinite multiplicity, that is,
    $$
A(\G) = \{ \varphi \in L^\infty(\G) : \varphi = (\xi,\eta)_{\lambda^\infty} \}.
    $$
    Moreover, $A(\G)$ is a closed $*$-ideal of $B(\G)$.
\end{prop}
We give some examples present in the work of Renault.
\begin{example}[{\cite[Example 1.2]{renault1997fourier}}]\label{Oid2group}
    When $\G = \Gamma$ is a discrete group, then $A(\G)$ reduces to the usual group Fourier algebra. 
\end{example}
\begin{example}[{\cite[Example 1.1]{renault1997fourier}}]
    When $\G = X $ is a groupoid reduced to its unit unit space then $A(\G,\mu) = B(\G,\mu) = L^\infty(X,\mu)$. 
\end{example}
\begin{example}[{\cite[Example 1.3]{renault1997fourier}}]
    Let $X$ be a locally compact space endowed with a probability measure $\mu$ and $\G = X \times X$ the trivial equivalence relation on $X$. Then $A(\G) = B(\G)$ are the space of \textit{Hilbertian functions} as defined by A. Grothendieck in \cite{grothendieck1956resume}, that is, the space of functions of the form $\varphi(x,y) = \langle \xi(x),\eta(y) \rangle_{L^2(X)}$ where $x,y \in X$ and $\xi,\eta: X \rightarrow L^2(X)$.
\end{example}
The following definition will be useful to prove that certain functions belong in $A(\G)$.
\begin{defi}
    A function $\varphi \in L^\infty(\G)$ is said to have \textit{r-compact support} if for every compact $K \subseteq X$, the set $\supp(\varphi) \cap K$ is compact.
\end{defi}
\begin{prop}[{\cite[Proposition 1.6]{renault1997fourier}}]\label{rcompactsupplemma}
 Let $\varphi \in B(\G)$, if $\varphi$ has r-compact support, then $\varphi \in A(\G)$. Moreover, the set of such functions is dense in $A(\G)$.
\end{prop}
\section{Duality and Multipliers}
We now discuss multipliers of the Fourier algebra and relate them to maps on the operator algebras associated with $\G$. As before, we fix a Hausdorff étale measured groupoid $(\G,\mu)$ with unit space $X := \gob$.
\begin{defi}\label{Groupoid Multipliers Definiiotn}
    Let $\varphi \in L^\infty(\G)$, $\varphi$ is said to be a \textit{multiplier of the Fourier algebra} if the map $M_\varphi: A(\G) \rightarrow A(\G)$ given by 
    $$
M_\varphi\psi = \varphi\psi
    $$
    is well defined and bounded. The set of such $\varphi$ is denoted by $M\!A(\G)$. Moreover, if $M_\varphi$ is completely bounded, $\varphi$ is said to be a \textit{completely bounded multiplier} and we denote the set of such multipliers by $M_0A(\G)$ which we endow with the norm
    $$
\|\varphi\|_{M_0A(\G)} = \|M_\varphi\|_{\text{cb}}.
    $$
\end{defi}
Furthermore, if $\varphi \in A(\G)$ the multiplication map $\psi \mapsto \varphi\psi$ is obviously completely bounded with norm less than $\|\varphi\|_{B(\G)}$ and this generalizes to $B(\G)$ giving the following proposition.
\begin{prop}[{\cite[Propositition 3.3]{renault1997fourier}}]\label{algebras inclusion}
We have the following norm-decreasing inclusions:
$$
A(\G) \subseteq B(\G) \subseteq M_0A(\G) \subseteq M\!A(\G).
$$
\end{prop}
As with the group case, $B(\G),A(\G)$ have a duality theory, which is more involved and we will present it but not go in detail. We start by observing that $C^*(\G,\mu),C_r^*(\G,\mu),L^2(X),L^2(X)^*$ are naturally operator spaces, where $L^2(X)^*$ is the dual operator space. Moreover, they are operator modules over $L^\infty(X)$ with actions defined similarly to what we did before:
\begin{align*}
    h * f * k(\gamma) := h(r(\gamma))f(\gamma)k(d(\gamma)) \quad &\text{for} \quad h,k \in L^\infty(X),f \in C_c(\G), \\
    ha(x) = h(x)a(x) \quad \quad &\text{for} \quad h \in L^\infty(X),a \in L^2(X), \\
    a^*h = (\overline{h}a)^* \quad \quad &\text{for} \quad h \in L^\infty(X),a \in L^2(X).
\end{align*}
See \cite[p.464]{renault1997fourier} for a more elaborate discussion. We define the spaces 
$$
\mathcal{X}(\G,\mu) := L^2(X)^* \mhotimes C^*(\G,\mu) \mhotimes L^2(X)
$$
and 
$$
\mathcal{X}_r(\G,\mu) := L^2(X)^* \mhotimes C_r^*(\G,\mu) \mhotimes L^2(X)
$$
where $\mhotimes$ is the module Haagerup tensor product over $L^\infty(X)$ (see \cite{blecher2000categories} for the definition of this tensor product). We will write $\mathcal{X}(\G)$ and $\mathcal{X}_r(\G)$ throughout the reminder of this section.
\begin{prop}[{\cite[Proposition 2.1]{renault1997fourier}}]\label{representX(G)}
    Let $u \in \mathcal{X}(\G)$ be an arbitrary element. Then there exist $a,b \in L^2(X)$ and $f \in C^*(\G,\mu)$ such that 
    $u = a^* \otimes f \otimes b$. Moreover, 
    $$
\|u\| = \inf \{ \|a\|_2\|f\|_\mu\|b\|_2: u= a^* \otimes f \otimes b \text{ where }a,b \in L^2(X) \text{ and } f \in C^*(\G,\mu) \}.
    $$
\end{prop}
If $\lambda: C^*(\G,\mu) \rightarrow C_r^*(\G,\mu)$ is the canonical quotient map, we can define $\Lambda: \id \otimes \lambda \otimes \id: \mathcal{X}(\G) \rightarrow \mathcal{X}_r(\G)$ which will also be a quotient in the sense of operator spaces. In particular, we obtain a reduced analog of Proposition \ref{representX(G)}.
\begin{prop}\label{representXr(G)}
Let $u \in \mathcal{X}_r(\G)$ be an arbitrary element. Then there exist $a,b \in L^2(X)$ and $f \in C^*_r(\G,\mu)$ such that 
    $u = a^* \otimes f \otimes b$. Moreover, 
    $$
\|u\| = \inf \{ \|a\|_2\|f\|_r\|b\|_2: u= a^* \otimes f \otimes b \text{ where }a,b \in L^2(X) \text{ and } f \in C^*_r(\G,\mu) \}.
    $$    
\end{prop}
We now see that elements of $B(\G)$ are related to functionals on $\mathcal X(\G)$.
\begin{prop}[{\cite[Lemma 2.1]{renault1997fourier}}]\label{groupoiddualityB(G)lemma}
    Let $\varphi = (\xi,\eta)_L \in B(\G)$, where $(\hil,L)$ is a $\G$-Hilbert bundle and $\xi,\eta \in L^\infty(X,\hil)$. Then there exists a unique bounded linear functional $\Phi$ on $\mathcal X(\G)$ such that for $a,b \in L^2(X)$ and $ f\in C_c(\G)$,
    $$
\Phi(a^* f b) = \int_{\G} \overline{a} * \varphi f * b \, d\nu_0 = \langle a \cdot \xi,\pi_L(f)(b \cdot \eta) \rangle_{L^2(\G,\nu)}, 
    $$
    where $a \cdot \xi(x) := a(x)\xi(x)$.
\end{prop}
\begin{theorem}[{\cite[Theorem 2.1]{renault1997fourier}}]\label{theoremdualB(G)}
    The correspondence $ \varphi \mapsto \Phi$ of Proposition \ref{groupoiddualityB(G)lemma} is an isometric isomorphism $\Xi: B(\G) \rightarrow \mathcal X(\G)^*$.
\end{theorem}
We will prove the above map induces an isometry $A(\G) \rightarrow \mathcal{X}_r(\G)^*$.
\begin{corollary}\label{embeddA(G)}
    Let $\varphi \in A(\G)$ and $\Phi = \Xi(\varphi) \in \mathcal{X}(\G)^*$ the associated functional of Theorem \ref{theoremdualB(G)}, then $\Phi$ factors through $\Lambda:\mathcal{X}(\G) \rightarrow \mathcal{X}_r(\G) $. Thus, we obtain a functional $\Phi_r \in \mathcal{X}_r(\G)^*$ and the correspondence $\Xi_r : A(\G) \ni \varphi \mapsto \Phi_r \in \mathcal{X}_r^*(\G)$ is an isometry.
\end{corollary}
\begin{proof}
Assume that $\varphi \in \A(\G)$ and let $\xi,\eta \in L^\infty(X,\Li^2(\G))$ be such that $\varphi = (\xi,\eta)_\lambda$. If $f \in C_c(\G)$ and $a,b \in L^2(X)$, then by Proposition \ref{groupoiddualityB(G)lemma}
\begin{align*}
    \Phi(a^*fb) = \langle a \cdot \xi,\lambda(f)(b \cdot \eta) \rangle_{L^2(\G,\nu)} \leq \|a\|_2\|\xi\|_\infty\|f\|_r\|\eta\|_\infty \|b\|_2.
\end{align*}
Thus we get a well defined functional $\Phi_r \in \mathcal{X}_r(\G)^*$ such that $\|\Phi_r\| = \|\Phi\| = \|\varphi\|_{B(\G)}$.
\end{proof}
This is not the end of the duality theory of $A(\G)$, we can go further. It is known that if $\G = \Gamma$ is a discrete group, $A(\Gamma)^* \cong \vN(\Gamma)$ thus $A(\Gamma)$ is the predual of $\vN(\Gamma)$. To generalize this, Renault introduced an operator space structure on $B(\G),A(\G)$ 
using the map $\Xi: B(\G) \rightarrow \mathcal{X}(\G)^*$ of Theorem \ref{theoremdualB(G)} and noting that there exists a complete isometry \cite[Proposition 2.2]{renault1997fourier}
$$
\mathcal{X}(\G)^* \cong \text{CB}_{X,X}(C^*(\G,\mu),B(L^2(X)))
$$
onto the space of completely bounded maps $C^*(\G,\mu) \rightarrow B(L^2(X))$ that commute with the left and right action of $L^\infty(X)$. Renault then defined matrix norm structures by letting 
$$
M_n(B(\G)) := \text{CB}_{X,X}(C^*(\G,\mu),B(L^2(X)) \otimes M_n)
$$
which turns $\Xi: B(\G) \rightarrow \mathcal{X}(\G)^*$ into a completely isometry. We have $M_n(A(\G)) \subseteq M_n(B(\G))$ and it can be seen that $\Xi_r : A(\G) \rightarrow \mathcal{X}_r(\G)^*$ is also a complete isometry. To compete the duality theory of $A(\G)$ in paralel to Eymard's result for gorups, Renault introduced the space
$$
\mathcal{Y}(\G) := L^2(X)^* \mhotimes A(\G) \mhotimes L^2(X)
$$
and proved it is the predual of $\vN(\G,\mu)$, see \cite[Theorem 2.3]{renault1997fourier}. In similarity with the group case, the main consequence of this duality theory is that if we have a multiplier $\varphi \in M\!A(\G)$, the multiplication map $M_\varphi:A(\G) \rightarrow A(\G)$ induces $1 \otimes M_\varphi \otimes 1: \mathcal{Y}(\G) \rightarrow \mathcal{Y}(\G)$ and by transposition we obtain $\overline{m}_\varphi := (1 \otimes M_\varphi \otimes 1)^*  :\vN(\G,\mu) \rightarrow \vN(\G,\mu)$ which is also defined through pointwise multiplication by $\varphi$. 

\begin{prop}\label{multipliers}
    Let $\varphi \in C_b(\G)$ (continuous and bounded functions), the following conditions are equivalent:
    \begin{enumerate}[(i)]
    \itemsep0pt
        \item Pointwise multiplication by $\varphi$ defines a (completely) bounded linear map from $A(\G)$ into itself with norm less than one; 
        \item Pointwise multiplication by $\varphi$ defines a (completely) bounded linear map from $\vN(\G)$ into itself with norm less than one.
        \item Pointwise multiplication by $\varphi$ defines a (completely) bounded linear map from $C_r^*(\G)$ into itself with norm less than one.
    \end{enumerate}
\end{prop}
\begin{proof}
    The equivalence of $(i) \iff (ii)$ is done in \cite[Propositions 3.1 and 3.2]{renault1997fourier} using the duality results mentioned above. The implication $(ii) \implies (iii)$ follows from the fact that $C_r^*(\G)$ is invariant for $m_\varphi$, meaning $m|_{C_r^*(\G)}: C_r^*(\G) \rightarrow C_r^*(\G)$ is a well defined operator that satisfies the norm requirements. We now prove $(iii) \implies (i) $, suppose $m_\varphi : C_r^*(\G) \rightarrow C_r^*(\G) $ is well defined and (completely) bounded with (c.b.) norm less than one and we let
    $$
    \mathcal{M}_\varphi := \id \otimes m_\varphi \otimes \id: \mathcal{X}_r(\G) \rightarrow \mathcal{X}_r(\G)
    $$
    then $\mathcal{M}_\varphi$ is bounded with the same (c.b.) norm. Dualizing we obtain $ \mathcal{M}_\varphi^*: \mathcal{X}_r(\G)^* \rightarrow \mathcal{X}_r(\G)^*$ which again is bounded with (c.b.) norm less than one. We define $S: A(\G) \rightarrow B(\G)$ to be the composition
    $$
A(\G) \xrightarrow[\quad \quad]{\textstyle{\Xi_r}} \mathcal{X}_r(\G)^* \xrightarrow[\quad \quad] {\mathcal{M}_\varphi^*}\mathcal{X}_r(\G)^* \xrightarrow[\quad \quad]{\textstyle{\Lambda^*}} \mathcal{X}(\G)^* \xrightarrow[\quad \quad]{\textstyle{\Xi^{-1}}} B(\G).
    $$
We claim that for $\psi \in A(\G)$ with $r$-compact support, we have $S(\psi) = \varphi\psi \in B(\G)$ which also has $r$-compact support, thus by Proposition \ref{rcompactsupplemma} $\varphi\psi \in A(\G)$ and $S$ is actually a map $A(\G) \rightarrow A(\G)$ such that $S = M_\varphi$. Since $S$ is defined as the composition of c.b. maps with (c.b.) norm less than one, we will conclude that $M_\varphi$ satisfies the norm requirements in (i). We now prove the claim, let $\psi \in A(\G)$ with $r$-compact support and $u = a^* \otimes f \otimes b \in \mathcal{X}(\G)$ where $a,b \in L^2(X)$ and $f \in C_c(\G)$. We define $\Psi := \Lambda^*\mathcal{M}_\varphi^* \, \Xi_r(\psi) $ and compute 
\begin{align*}
   \Psi(a^* \otimes f \otimes b) = \mathcal{M}_\varphi^* \, \Xi_r(\psi)(a^* \otimes \lambda(f) \otimes b) = \Xi_r(\psi)(a^* \otimes \varphi\lambda(f) \otimes b) 
\end{align*}
Since $f \in C_c(\G)$, then $\lambda(f) = f$ and $\varphi f \in C_c(\G)$, thus by Proposition \ref{representXr(G)}
$$
\Xi_r(\psi)(a^* \otimes \varphi\lambda(f) \otimes b)  = \int_{\G} \overline{a} * \psi \varphi f * b \, d\nu_0.
$$
By Theorem \ref{theoremdualB(G)}, there exists a $\phi \in B(\G)$ such that $\Xi(\phi) = \Psi$, meaning that for all $a,b \in L^2(X)$ and $f \in C_c(\G)$ we have
$$
\int_{\G} \overline{a} * \phi f * b \, d\nu_0 = \int_{\G} \overline{a} * \psi \varphi f * b \, d\nu_0.
$$
from which can deduce $\phi = \varphi \psi$ $\nu$-a.e. ($\nu$ and $\nu_0$ have the same null sets) and so $\varphi \psi \in B(\G)$. To conclude note that 
$$
\Lambda^*\mathcal{M}_\varphi^* \, \Xi_r(\psi) = \Xi(\varphi\psi) \iff S(\psi) = \Xi^{-1}\Lambda^*\mathcal{M}_\varphi^* \, \Xi_r(\psi) = \varphi\psi \in B(\G),
$$
as we wanted to show.
\end{proof}
\section{Weak Amenability}
Before defining weak amenability let us set the stage and compile the different notions of amenability for (measured) groupoids. For every type of amenability, there are many equivalent characterizations but we will only give one for each, according to our needs. For a full review of the subject, we refer to \cite{anantharaman2001amenable} and we will use \cite{anantharaman2023amenability} as a reference. We start with the measured groupoid setting.
\begin{defi}(Weak Godement Condition)
    A measured groupoid $(\G,\mu)$ is said to be \textit{amenable} if there exists a net $(\xi_i)$ of Borel functions on $\G$ such that for all $\gamma \in \G$, 
    $$
    \sum\limits_{\beta \in G^{r(\gamma)}} |\xi_i(\beta)|^2 = 1, \quad \text{for all } \gamma \in \G
    $$
    and the functions
    $$
    g_i(\gamma) := \sum\limits_{\beta \in \G^r(\gamma)} \overline{\xi_i(\beta)}\xi_i(\gamma^{-1}\beta)
    $$
    converge to $1$ in the weak-$*$ topology of $L^\infty(\G)$.
\end{defi}
The first condition tells us that the functions $\xi_i \in L^\infty(\gob, \Li_\mu^2(\G))$ with $\|\xi\|_\infty = 1$. It can be seen that $g_i := (\xi_i,\xi_i)_\lambda$, meaning $g_i \in A(\G)$ and $\|g_i\|_{B(\G)} \leq 1$. We now proceed to amenability in the topological setting.
\begin{defi}
    An étale groupoid is said to be \textit{measurewise amenable} if for all quasi-invariant measures $\mu$, $(\G,\mu)$ is amenable. Moreover, it is said to be \textit{topologically amenable} if there exists a net $(g_i) \in C_c(\G)$ of positive definite functions such that $\| g|_{\gob}\|_\infty \leq 1$ and $g_i \rightarrow 1$ uniformly on compact subsets.
\end{defi}
It is easy to see that topological amenability implies measurewise amenability and for étale groupoids, the converse also holds \cite{anantharaman2001amenable}. We now move on to weak amenability where we will also start from the measured groupoid case.
\begin{defi}
    Let $(\G,\mu)$ be a measured groupoid, then $(\G,\mu)$ is said to be \textit{weakly amenable} if there exists a net $(\varphi_i) \in A(\G)$ and a $C>0$ such that $\sup_i \|\varphi_i \|_{M_0A(\G)} = C < \infty$ and $\varphi_i \xrightarrow{} 1$ in the weak-$*$ topology of $L^\infty(\G)$. The \textit{Cowling-Haagerup constant} is the infimum of all $C > 0$ for which there is a net as above and is denoted by $\Lambda_{\text{cb}}(\G,\mu)$.
\end{defi}
When $(\G,\mu)$ is not weakly amenable, we set $\cb(\G,\mu) = + \infty.$
\begin{defi}
    Let $\G$ be an étale groupoid, then $\G$ is said to be \textit{measurewise weakly amenable} if there exists a $C > 0$ such that for all quasi-invariant probability measures $\mu$ on $\gob$, $(\G,\mu)$ is weakly amenable with $\cb(\G,\mu) \leq C$. The \textit{measurewise Cowling-Haagerup constant} is the infimum of all $C > 0$ for which above holds and is denoted by $\Lambda_{\text{mcb}}(\G)$. 
\end{defi}
Again, when $\G$ is not measurewise weakly amenable we set $\Lambda_{\text{mcb}}(\G) = + \infty$. By definition, we also have
$$
\Lambda_{\text{mcb}}(\G) := \sup \{ \cb(\G,\mu): \mu \text{ is a quasi-invariant probability measure on } \gob \}.
$$
For topological weak amenability, we will want to make use of continuous functions and since $\G$ is étale, we have that $C_c(\G) \subseteq A(\G)$, where we endow $\G$ with a quasi-invariant measure with full support.
\begin{defi}
    Let $\G$ be an étale groupoid, then $\G$ is said to be \textit{(topologically) weakly amenable} if there exists  a net $(\varphi_i) \in C_c(\G)$, a $C > 0$ and quasi-invariant probability measure with full support $\mu$ on $\gob$, such that $\sup_i \|\varphi_i \|_{M_0A(\G)} = C < \infty$ and $\varphi_i \xrightarrow{} 1$ uniformly on compact subsets. 
\end{defi}
Before defining the appropriate notion of the Cowling-Haagerup constant, let us see that it will not depend on the choice of $\mu$.
\begin{prop}\label{inequalitymeasures}
    Let $\mu_1,\mu_2$ be quasi-invariant measures on $\gob$ with $\supp(\mu_1) = \supp(\mu_2)$. Then for $\varphi \in C_b(\G)$, $\|\varphi\|_{M_0A(\G,\mu_1)} \leq 1$ if and only if $\|\varphi\|_{M_0A(\G,\mu_2)} \leq 1$. 
\end{prop}
\begin{proof}
    Since $\supp(\mu_1) = \supp(\mu_2)$ we have that the identity $C_c(\G) \rightarrow C_c(\G)$ extends to an isometric isomorphism $C_r^*(\G,\mu_1) \cong C_r^*(\G,\mu_2)$. The result now follows from applying Proposition \ref{multipliers}.
\end{proof}
\begin{defi}
    The \textit{(topological) Cowling-Haagerup constant} of $\G$ is the infimum of all $C > 0$ such that there exists a quasi-invariant measure $\mu$ and a net $(\varphi_i) \in C_c(\G)$ such that $\sup_i \|\varphi_i \|_{M_0A(\G)} = C < \infty$ and $\varphi_i \xrightarrow{} 1$ uniformly on compact subsets. It is denoted by $\cb(\G)$. 
\end{defi}
Of course, if $\G$ is not topologically weakly amenable we set $\cb(\G) = + \infty$.
\begin{remark}
    This work will focus on topological weak amenability so from here on the term "weak amenability" will mean this.
\end{remark}
From the definitions, it is easily seen that amenable groupoids are weakly amenable with the associated Cowling-Haagerup constants equal to one. We state this more precisely in the following proposition.
\begin{prop}
    Let $\G$ be an étale groupoid. Then:
    \begin{enumerate}
        \itemsep-5pt
        \item If $\mu$ is a quasi-invariant measure and $(\G,\mu)$ is amenable, then $(\G,\mu)$ is weakly amenable and $\cb(\G,\mu) = 1$;
        \item If $\G$ is measurewise amenable, then $(\G,\mu)$ is measurewise weakly amenable and\\
        $\Lambda_{\text{mcb}}(\G) = 1$;
        \item If $\G$ is amenable, then $\G$ is weakly amenable and $\cb(\G) = 1$.
    \end{enumerate}
\end{prop}
On the other hand, there is very relevant necessary condition for topological weakly amenable groupoids which we describe now. Let $F \subseteq X$ be a closed invariant subset and $U:= X \backslash F$. Recall from section two that there exists a sequence
\begin{equation}\label{exactredreduction}
0 \longrightarrow C^*_r(\G|_U) \xrightarrow[\quad\,\,]{\iota} C^*_r(\G) \xrightarrow[\quad\,\,]{p} C^*_r(\G|_F) \longrightarrow 0,    
\end{equation}
which is not necessarily exact (in the middle). However, if the sequence in Equation \ref{exactredreduction} is exact for all closed invariant subsets $F \subseteq X$, then $\G$ is said to be \textit{inner exact}. It is well known that there exists a norm-decreasing injective map $j: C_r^*(\G) \rightarrow C_0(\G)$ such that for all $f \in C_c(\G)$, $j(f) = f$ \cite[Proposition 3.3]{sims2017etale}. This allows us to define the \textit{support} of an element $f \in C_r^*(\G)$ by 
$$
\supp(f) := \supp(j(f)).
$$
Furthermore, we get a commutative diagram.
\[\begin{tikzcd}
	0 & {C_r^*(\G|_U)} & {C_r^*(\G)} & {C_r^*(\G|_F)} & 0 \\
	& {C_0(\G|_U)} & {C_0(\G)} & {C_0(G|_F)}
	\arrow[from=1-1, to=1-2]
	\arrow["{\textstyle{\iota}}", from=1-2, to=1-3]
	\arrow["{\textstyle{j_U}}"', from=1-2, to=2-2]
	\arrow["{\textstyle{p}}", from=1-3, to=1-4]
	\arrow["{\textstyle{j}}", from=1-3, to=2-3]
	\arrow[from=1-4, to=1-5]
	\arrow["{\textstyle{j_F}}", from=1-4, to=2-4]
	\arrow["{\textstyle{\iota'}}", from=2-2, to=2-3]
	\arrow["{\textstyle{p'}}", from=2-3, to=2-4]
\end{tikzcd}\]
The map $p':C_0(\G) \rightarrow C_0(\G|_F)$ is simply given by function restriction so 
$$
\ker(p') = \{ f \in C_0(\G): \supp(f) \subseteq \G \backslash \G|_F = \G|_U \}.
$$
Since $j,j_F$ is injective, we also get that
$$
\ker(p) = \{ f \in C_r^*(\G): \supp(f) \subseteq \G|_U \}.
$$
\begin{prop}\label{WA implies Inner Exact}
    Let $\G$ be a Hausdorff étale groupoid. If $\G$ is weakly amenable, then it is inner exact.
\end{prop}
\begin{proof}
Let $F$ be a closed invariant subset and $U:= X \backslash F$. We have to show that the sequence from Equation \ref{exactredreduction} is exact. Let $f \in C_c(\G|_U)$, then $p(\iota(f)) = 0$ and so $\im(\iota) \subseteq \ker(p)$. For the reverse inclusion, we will have to use the assumption that $\G$ is weakly amenable, thus let $(\varphi_i)_{i \in I} \in C_c(\G)$, $C > 0$ and $\mu$ a quasi invariant measure testifying this property for $\G$. Let $f \in \ker(p)$, then for all $i \in I$,
$$
\supp(\varphi_i f) \subseteq \supp (\varphi) \cap \supp(f)
$$
meaning $\varphi_i f \in C_c(\G|_U)$ and so $\varphi_i f \in \im(\iota)$. The fact that $\varphi_i f \rightarrow f$ follows from an argument done in the proof of Proposition \ref{weakinequality}.
\end{proof}
In his paper \cite{willett2015non}, Willet proves that there exist non-amenable groupoids that have the weak containment property. To do so, he used the class of HLS groupoids of Higson, Lafforgue and Skandalis \cite{higson2002counterexamples}. Willet proved that an HLS groupoid $\G$ is amenable if and only if it is inner exact. Thus, the previous proposition and Willet's results tell us that:
\begin{prop}
    There exist non weakly amenable groupoids that have weak containment. Conversely, it is known that there exist weakly amenable groups that do not have weak containment.
\end{prop}
Let $\G$ be a weakly amenable étale groupoid and $\mu$ be an arbitrary quasi-invariant measure on $X$, then $\Omega := \supp(\mu)$ is a closed invariant subset of $X$. By Proposition \ref{WA implies Inner Exact}, the sequence (where $U := X \backslash \Omega)$ 
\begin{equation}
0 \longrightarrow C^*_r(\G|_U) \xrightarrow[\quad\,\,]{\iota} C^*_r(\G) \xrightarrow[\quad\,\,]{p} C^*_r(\G|_\Omega) \longrightarrow 0,    
\end{equation}
is exact. Recalling that $C_r^*(\G|_\Omega) \cong C_r^*(\G,\mu)$ we conclude that $ C_r^*(\G,\mu)$ is a quotient of $C_r^*(\G)$ with kernel $C_r^*(\G|_U)$.
\begin{theorem}
    Let $\G$ be a weakly amenable étale groupoid, then $\G$ is measurewise weakly amenable and $\text{m}\cb(\G) \leq \cb(\G)$. 
\end{theorem}
\begin{proof}
    Let $\mu_0$ be a quasi-invariant measure with full support, $(\varphi_i) \in C_c(\G)$ a net and $C>0$ testify the weak amenability of $\G$, let also $\mu$ be an arbitrary quasi-invariant measure. We prove that for all $i \in I$, $\|\varphi_i\|_{M_0A(\G,\mu)} \leq C$. By Proposition \ref{multipliers}, multiplication by $\varphi_i$ induces an operator $m_{\varphi_i} : C_r^*(\G) \rightarrow C_r^*(\G)$. Let $p: C_r^*(\G) \rightarrow C_r^*(\G|_\Omega) \cong C_r^*(\G,\mu)$ be the extension of the restriction map, then $p(\varphi_i) \in C_c(\G|_\Omega)$ and so again by Proposition \ref{multipliers} we obtain a linear operator $m_{p(\varphi_i)}:C_r^*(\G,\mu) \rightarrow C_r^*(\G,\mu)$ and a commutative diagram 
    \[\begin{tikzcd}
	{C_r^*(\G)} & {C_r^*(\G|_\Omega)} \\
	{C_r^*(\G)} & {C_r^*(\G|_\Omega)}
	\arrow["{\textstyle{p}}", from=1-1, to=1-2]
	\arrow["{\textstyle{m_{\varphi_i}}}"', from=1-1, to=2-1]
	\arrow["{\textstyle{m_{p(\varphi_i)}}}", from=1-2, to=2-2]
	\arrow["\textstyle{p}"', from=2-1, to=2-2]
\end{tikzcd}\]
from which is possible to deduce that $\|m_{p(\varphi_i)}\|_{\text{cb}} \leq \|m_{\varphi_i}\|_{\text{cb}}$. Using Proposition \ref{multipliers} once more we conclude $\|\varphi_i\|_{M_0A(\G,\mu)} \leq C$. To finish, we note that uniform convergence in compact sets implies weak-$*$ convergence in $L^\infty(\G)$ and so $(\G,\mu)$ is weakly amenable with $\cb(\G,\mu) \leq C$. Since $\mu$ was arbitrary, we get that $\Lambda_{\text{mcb}} \leq C$. 
\end{proof}
\section{Completely Bounded Approximation Property of Cartan Pairs}
As always, we let $(\G,\mu)$ be a Hausdorff measured étale groupoid with unit space $X := \gob$. Before we move on to the approximation properties on the algebraic side, let's see what to expect. Suppose $\varphi \in C_c(\G)$ and consider the induced multiplier $m_\varphi : C_r^*(\G,\mu) \rightarrow C_r^*(\G,\mu)$, unlike what happens for the group case, it is not necessary that $m_\varphi$ has finite rank, indeed, we can just take $X$ to be a compact Hausdorff space, $\mu$ be a Borel probability measure on $X$ with full support and $\G = \gob = X$. The constant $1$ function is in $ A(X) = L^\infty(X,\mu)$ and will induce the identity map on $C_r^*(\G,\mu) = C(X)$ which does not have finite rank (unless $X$ is finite). To solve the above problem, we introduce the correct algebraic setting based on \cite{renault2024cartan}.
\begin{defi}
    Let $A \subseteq B$ be an abelian sub $C^*$-algebra of a $C^*$-algebra $B$. We will say that $A$ is a \textit{quasi Cartan subalgebra} if:
    \begin{enumerate}[(i)]
    \itemsep-5pt
        \item $A$ contains an approximate unit of $B$;
        \item $A$ is regular;
        \item there exists a conditional expectation $E:B \rightarrow A$. 
    \end{enumerate}
    Then $(B,A)$ is called a \textit{quasi Cartan pair}. If moreover 
    \begin{enumerate}[(iv)]
        \item $A$ is maximal abelian.
    \end{enumerate}
    Then $(B,A)$ is called a \textit{Cartan pair}.
\end{defi}
It is known that if $\G$ is an étale essentially principal groupoid in the sense of \cite[Definition II.4.3]{renault2006groupoid}, then $(C_r^*(\G),C_0(X))$ is Cartan pair \cite{renault2024cartan}. Moreover, if $\mu$ is a quasi-invariant measure on $X$ with support $\Omega := \supp(\mu)$ then $\G|_\Omega$ is essentially principal and $(C_r^*(\G|_\Omega),C_0(\Omega)) = (C_r^*(\G,\mu), C_0(\Omega))$ is a Cartan pair. The fact that $\G$ is assumed to be essentially principal is there to ensure that $C_0(F) \subseteq C_r^*(\G|_F)$ is maximal abelian for all closed invariant subsets $F \subseteq X$. However, we will not need this fact so we focus on quasi Cartan pairs, which by the following proposition will be enough.
\begin{prop}
    Let $(\G,\mu)$ be a Hausdorff étale measured groupoid and $\Omega := \supp(\mu)$. Then $(C_r^*(\G|_\Omega),C_0(\Omega))$ is a quasi Cartan pair. 
\end{prop}
\begin{proof}
    Since $\G$ is étale, then $\G|_\Omega$ also is and by \cite[Corollary 4.8]{renault2024cartan} $C_0(\Omega)$ is regular. Moreover, since $\G|_\Omega$ is étale Hausdorff, the restriction map $C_c(\G|_\Omega) \rightarrow C_c(\Omega)$ extends to a conditional expectation $E: C_r^*(\G|_\Omega) \rightarrow C_0(\Omega)$ \cite[p.204]{brown2008textrm}. With respect to condition (i), simply take positive compactly supported functions $e_i \in C_c(\Omega)$ such that $e_i \rightarrow \chi_\Omega$ uniformly on compact subsets of $\gob$.
\end{proof}
When $(B,A)$ is a quasi Cartan pair, $B$ becomes a (right) pre-Hilbert $C^*$-module over $A$ with inner product given by $\langle a,b \rangle_A = E(a^*b)$ for $a,b \in B$. The Hilbert $C^*$-module obtained after separation and completion is denoted by $L^2(B,E)$. The induced norm satisfies the following inequality for $b \in B$
$$
\|b\|_{L^2(B,E)}^2 = \|\langle b,b \rangle_A \| = \|E(b^*b)\| \leq \|b\|^2. 
$$
This, together with the other conditions of quasi Cartan subalgebras, ensures that there exists a norm decreasing inclusion $B \rightarrow L^2(B,E)$.
\begin{defi}
    Let $B$ be a right pre-Hilbert $C^*$-module over $A$ and be $T: B \rightarrow B$ be an $A$-linear operator. $T$ is said to have \textit{$A$-rank one} if there exist $b_1,b_2 \in B$ such that for all $x \in B$,
    $$
    Tx =  b_1 \cdot \langle b_2,x \rangle_A.
    $$
    When this is the case, we write $T = \Theta_{b_1,b_2}$. Furthermore, an $A$-linear map $T$ is said to have \textit{finite $A$-rank} if it is the sum of finitely many $A$-rank one maps.
\end{defi}
If $(B,A)$ is a quasi Cartan pair, there might be confusion on which norm we consider on $B$, but in this work we will focus on the norm already coming from $B$, meaning $\| \cdot \|_{L^2(B,E)}$ will not not be very relevant to us. With this in mind, when we consider $A$-linear bounded operators $T: B \rightarrow B$, they are bounded on the norm of $B$.
\begin{prop}\label{decomposablerankone}
    Let $(B,A)$ be a quasi Cartan pair, then for every $b_1,b_2 \in B$, the map $\Theta_{b_1,b_2}$ is completely bounded and $\|\Theta_{b_1,b_2}\|_{\text{cb}} \leq \|b_1\|\|b_2\|$.
\end{prop}
\begin{proof}
    For every $b \in B$, the map $m_b:B \rightarrow B$ given by left multiplication is c.b. and $\|m_b\|_{\text{cb}} \leq \|b\|$. Note that $\Theta_{b_1,b_2} =m_{b_1} \circ E \circ m_{b_2^*} $, meaning it is a composition of c.b. maps and we have
    $$
 \|\Theta_{b_1,b_2}\|_{\text{cb}} \leq \|b_1\|\|b_2\|.
    $$
\end{proof}
We now generalize the completely bounded approximation property to quasi Cartan pairs.
\begin{defi}
  A quasi Cartan pair $(B,A)$ has the \textit{completely bounded approximation property} (abbreviated by CBAP) if there exists a net $T_i: B \rightarrow B$ of $A$-linear operators with finite $A$-rank such that $T_i \xrightarrow{\text{SOT}} \id_B$ and $\sup_i \| T_i \|_{\text{cb}} = C < \infty$. The \textit{Cowling-Haagerup constant} is the least $C$ such that there is a net as above and is denoted by $\Lambda_{\text{cb}}(B,A)$. 
\end{defi}
Our goal now is to prove that multipliers of compactly supported functions have finite $C_0(X)$-rank.
\begin{prop}\label{finiterankmultipliers}
    Let $\varphi \in C_c(\G)$ and $m_\varphi:C_r^*(\G,\mu) \rightarrow C_r^*(\G,\mu)$ be the induced multiplier, then:
    \begin{enumerate}[(i)]
    \itemsep0pt
        \item If $\varphi$ is supported on a bisection $U \in \G^{op}$, then $m_\varphi$ has $C_0(X)$-rank one;
        \item $m_\varphi$ has finite $C_0(X)$-rank. 
    \end{enumerate}
\end{prop}
\begin{proof}
    For (i), let $\varphi \in C_c(U)$, where $U \in \G^{op}$. We will show that $m_\varphi = \Theta_{\varphi,\chi_U}$, let $f \in C_c(\G)$ and $x \in X$ and first observe that 
    $$
\langle \chi_U,f \rangle_A (x) = \sum\limits_{\gamma \in \G_x} \chi_U(\gamma)f(\gamma) = \begin{cases}
    f \circ d_U^{-1}(x) \quad &\text{ if } x \in d(U), \\
    0 &\text{ otherwise.}
\end{cases}
    $$
    so if $f \in C_c(\G)$ and $\gamma \in U$ we get 
    \begin{align*}
  \Theta_{\varphi,\chi_U}(f)(\gamma) = \varphi * \langle \chi_U,f \rangle_A (\gamma) = \varphi(\gamma) \,  \langle \chi_U,f \rangle_A (d(\gamma)) 
  = \varphi(\gamma) \,  f \circ d_A^{-1}(d(\gamma)) &= \varphi(\gamma)f(\gamma) \\
  &=m_\varphi(f)(\gamma).  
    \end{align*}
Of course if $\gamma \notin U$, then $\Theta_{\varphi,\chi_U}(f)(\gamma) = m_\varphi(f)(\gamma) = 0$ and we conclude that $\Theta_{\varphi,\chi_U} = m_\varphi$ on $C_c(\G,\mu)$ and hence on all of $C_r^*(\G)$, proving (i). Point (ii) follows from a partition of unity argument, if $\varphi \in C_c(\G)$, we cover $\supp(\varphi)$ with open bisections $U_1,...,U_N \in \G^{op}$ and take $\psi_1,...,\psi_N$ to be a partition of unity subordinate to $U_1,...,U_N$. Letting $\varphi_i := \varphi \psi_i \in C_c(U_i) \subseteq C_c(\G)$ for $i = 1,...,N$, we have that 
$$
m_{\varphi} = \sum\limits_{i=1}^N m_{\varphi_i},
$$
by (i) each $m_{\varphi_i}$ is an $A$-rank one map and so we conclude that $m_{\varphi}$ has finite $C_0(X)$-rank.
\end{proof}
\begin{prop}\label{weakinequality}
    Let $\G$ be an étale groupoid, then 
    $$
\cb(C_r^*(\G),C_0(\gob)) \leq \cb(\G).
    $$
\end{prop}
\begin{proof}
    If $\G$ is not weakly amenable, there is nothing to prove. Now let $\mu$ be a quasi-invariant measure on $X$ with full support, $(\varphi_i) \in C_c(\G)$ be a net and $C > 0$ testify the (topological) weak amenability of $\G$. By Propositions \ref{multipliers} and $\ref{finiterankmultipliers}$ the operators $T_i := m_{\varphi_i}: C_r^*(\G) \rightarrow C_r^*(\G)$ have finite $C_0(X)$-rank and satisfy $\sup_i \|T_i\|_{\text{cb}} \leq C$, thus it only remains to show that $T_i \xrightarrow{\text{SOT}} \id_{C_r^*(\G)}$. Let $f \in C_c(\G)$ and cover the support of $f$ with finitely many open bisections $U_1,...,U_N \in \G^{op}$. Let $\psi_1,...,\psi_N$ be a partition of unity subordinate to $U_1,...,U_N$ and let $g_k := \psi_kf$ then
    \begin{align*}
    \| T_i(f) - f\|_r = \left\| \sum\limits_{k = 1}^N \varphi_ig_k - g_k \right\|_r \leq \sum\limits_{k = 1}^N \|\varphi_ig_k - g_k\|_r   
    \end{align*}
    Since $\varphi_ig_k - g_k \in C_c(U_k)$ is supported on a bisection and $\varphi_i \rightarrow 1$ uniformly on compact sets we conclude 
    $$
\sum\limits_{k = 1}^N \|\varphi_ig_k - g_k\|_r = \sum\limits_{k = 1}^N \|\varphi_ig_k - g_k\|_\infty \rightarrow 0.
    $$
    Thus $T_i \xrightarrow{\text{SOT}} \id_{C_r^*(\G)}$ and so $\cb(C_r^*(\G),C_0(X)) \leq \cb(\G)$.
\end{proof}
I was not able to prove the converse in general, however the next two sections are dedicated to classes of groupoids where I achieved the reverse inequality.
\section{Weak Amenability of Partial Actions}
In the previous section, we proved that for a Hausdorff étale groupoid $\G$ with unit space $X := \gob$, we have $\cb(C_r^*(\G),C_0(X)) \leq \cb(\G)$. This section is about the first class of groupoids to which we could prove the converse. More concretely, we will show that if $\G$ is a groupoid arising from a partial action of a discrete group $\Gamma$ on a locally compact space $X$ then 
\begin{equation}\label{equalitypartialactionintro}
\cb(\G) = \cb(C_r^*(\G),C_0(X)).    
\end{equation}
In Section 7.1, we define $C^*$-algebraic partial actions. As it is with regular group actions, it gives rise to a crossed product $C^*$-algebra $A \rtimes_r \Gamma$. We state Fell's absorption associated with these actions and finish the section with the existence of an injective $*$-homomorphism $A \rtimes_r \Gamma \xrightarrow{} C_r^*(\Gamma) \otimes A \rtimes_r \Gamma$. In Section 7.2, we turn to topological partial actions and the groupoid induced and in Section 7.3 we prove the equality in Equation \eqref{equalitypartialactionintro}.
    \subsection{Partial Actions}
We start with partial actions in the $C^*$-algebra setting.
\begin{defi}
    A ($C^*$-algebraic) \textit{partial action} of $\Gamma$ on a $C^*$-algebra $A$ is a pair
    $$
\beta = ( \{A_t\}_{t \in \Gamma}, \{\beta_t\}_{t \in \Gamma} )
    $$
    where each $A_t$ is a closed two-sided ideal of $A$ and each $\beta_t: A_{t^{-1}} \rightarrow A_t$ is a $*$-isomorphism. A \textit{($C^*$-algebraic) partial dynamical system} will mean the triple $(A,\Gamma,\beta)$.
\end{defi}
This information can also be coded in the form of a Fell bundle as it is done in \cite{exel2017partial}.
Its total space is 
$$
\A_\beta := \{(a,t) \in A \times \Gamma: t \in \Gamma \text{ and } a \in A_t \}.
$$
We will denote the element $(a,t)$ by $a\delta_t$, thus the fibers of the bundle can be written as 
$$
A_t = \{a\delta_t : a \in A_t \}.
$$
We have the skewed product
$$
(a\delta_s)(b\delta_t) = \beta_s(\beta_{s^{-1}}(a)b)\delta_{st}
$$
and involution
$$
(a\delta_s)^* = \beta_{s^{-1}}(a^*)\delta_{g^{-1}},
$$
for all $s,t \in \Gamma$ and $a \in A_s,b \in A_t$. With these operations, $\A_\beta$ becomes a Fell bundle, called the \textit{semi-direct product bundle relative to $\beta$} \cite[Proposition 16.6]{exel2017partial}. Denoting the space of compactly supported sections by $C_c(\A_\beta)$, it is known that $A$ embedds in $C_c(\A_\beta)$, so we can introduce an $A$-valued inner product defined by 
$$
\langle f,g \rangle = \sum\limits_{t \in \Gamma} f(t)^*g(t), \quad \text{ for all } f,g \in C_c(\A_\beta). 
$$
The completion of $C_c(\A_\beta)$ under the norm $\|\cdot\|_2$ induced by the inner product above is denoted by $\ell^2(\A_\beta)$. Moreover, for each $t \in \Gamma$ and $a \in A_t$ we define $ \lambda_t(a): \ell^2(\A_\beta) \rightarrow \ell^2(\A_\beta)$ by 
$$
C_c(\A_\beta) \ni (b\delta_u) \mapsto \beta_t(\beta_{t^{-1}}(a)b)\delta_{tu} \in C_c(\A_\beta)
$$
it is a bounded and adjointable operator \cite[Proposition 17.3]{exel2017partial} and the collection $\{\lambda_t\}_{t \in \Gamma}$ forms a representation of $\A_\beta$ on $B(\ell^2(\A_\beta))$ \cite[Proposition 17.4]{exel2017partial}, called the \textit{left regular representation of $\A_\beta$}. This induces a faithful representation on the space of sections of the bundle $\Lambda : C_c(\A_\beta) \rightarrow B(\ell^2(\A_\beta))$ \cite[Proposition 17.9]{exel2017partial}. The \textit{reduced crossed product $C^*$-algebra} $A \rtimes_r \Gamma$ is defined to be the closure of the image of $\Lambda$ in $B(\ell^2(\A_\beta))$, that is,
$$
A \rtimes_r \Gamma := \overline{\Lambda(C_c(\A_\beta))} \subseteq B(\ell^2(\A_\beta)).
$$
Now suppose $\pi = \{\pi_t\}_{t \in \Gamma}$ is a representation of $\A_\beta$ in a Hilbert space $H$ and let also $\lambda^\Gamma: \Gamma \rightarrow B(\ell^2(\Gamma))$ be the left regular representation of $\Gamma$. We can define $ \lambda^\Gamma \otimes \pi := \{ \lambda^\Gamma_t \otimes \pi_t\}_{t \in \Gamma},$ which will be a new representation of $\A_\beta$, this time on $\ell^2(\Gamma) \otimes H$. Denoting the induced representation of $C_c(\A_\beta)$ also by $\lambda^\Gamma \otimes \pi$ we obtain a Fell's absorption principle for Fell bundles\footnote{Although we are restricting to the case where the Fell bundle is given by a partial action, this is valid for a general one.}.
\begin{prop}[{\cite[Proposition 18.4]{exel2017partial}}]
Let $\pi$ be a representation of $\A_\beta$ on a Hilbert space $H$. Let $\lambda^\Gamma \otimes \pi: C_c(\A_\beta) \rightarrow B(\ell^2(\Gamma) \otimes H)$ be the induced representation as above, then $\lambda^\Gamma \otimes \pi$ vanishes on the kernel of $\Lambda$ meaning it descends to a representation $\lambda^\Gamma \otimes \pi: A \rtimes_r \Gamma \rightarrow B(\ell^2(\Gamma) \otimes H)$. Moreover, if $\pi_1$ is faithful, then so is $\lambda^\Gamma \otimes \pi: A \rtimes_r \Gamma \rightarrow B(\ell^2(\Gamma) \otimes H)$.
\end{prop}
The major consequence we will make use in this work is the existence of the following diagonal coaction 
\begin{corollary}\label{emb Fell version}
    There exists an injective $*$-homomorphism $A \rtimes_r \Gamma \rightarrow C_r^*(\Gamma) \otimes A \rtimes_r \Gamma$ such that for all $t \in \Gamma$ and $a \in A_t$, we have
    $$
\lambda_t(a) \mapsto \lambda^\Gamma_t \otimes \lambda_t(a).
    $$
\end{corollary}

    \subsection{Partial Action Groupoids}
We now give the definition of a topological partial action, which will be the ones that have a groupoid associated.
\begin{defi}
    A \textit{partial action} of $\Gamma$ on $X$ is a pair 
    $$
\alpha = (\{D_t\}_{t \in \Gamma}, \{\alpha_t\}_{t \in \Gamma}) 
    $$
    consisting of a collection $\{D_t\}_{t \in \Gamma}$ of open subsets of $X$ and a collection of homeomorphisms 
    $
\alpha_t : D_{t^{-1}} \rightarrow D_t,
    $
    such that 
    \begin{enumerate}[(i)]
    \itemsep0pt
        \item $D_1 = X$, and $\alpha_1$ is the identity map,
        \item $\alpha_s \circ \alpha_t \subseteq \alpha_{st}$, for all $g$ and $h$ in $\Gamma$.
    \end{enumerate}
\end{defi}
As it is the case with global actions, a partial action has an associated groupoid denoted by $X \rtimes \Gamma$, whose construction plays out very similarly. Its underlying set is 
$$
X \rtimes \Gamma := \{ (x,t): x \in D_t \} = \bigsqcup\limits_{t \in \Gamma} D_t \times \{t\},
$$
the source and range maps are $d(x,t) := (x,e)$ and $r(x,t) := (\alpha_t(x),e)$, respectively. Multiplication, when defined, is given by 
$$
(\alpha_s(x),t)(x,t) = (x,gh), \quad \text{ for all } s,t \in \Gamma \text{ and } x \in D_t.
$$
We give $X \rtimes \Gamma$ the subspace topology of $X \times \Gamma$ and since $\Gamma$ is discrete, $X \rtimes \Gamma$ will be étale. We moreover endow $X$ with a quasi-invariant measure $\mu$, obtaining a measured groupoid $(X \rtimes \Gamma, \mu)$. 
The (topological) partial action $\alpha$ of $\Gamma$ on $X$ induces a $C^*$-algebraic one in the usual way. We let $A_t := C_0(D_t)$ and $\beta_t: A_{t^{-1}} \rightarrow A_t$ be defined by 
$$
f \mapsto \beta_t(f)(x) = f(\alpha_{t^{-1}}(x)), \quad \text{ for all } x \in D_t.
$$
Define $\Delta: C_c(\A_\beta) \rightarrow C_c(X \rtimes \Gamma)$ by
$$
f \mapsto \Delta(f)(x,t) = f(t)(x), \quad \text{ for all } f \in C_c(\A_\beta), t\in \Gamma \text{ and } x \in D_t.
$$
\begin{prop}
    The map $\Delta: C_c(\A_\beta) \rightarrow C_c(X \rtimes \Gamma)$ is a $*$-isomorphism such that:
    \begin{enumerate}[(i)]
    \itemsep0pt
        \item It extends to a unitary operator $\ell^2(\A_\beta) \rightarrow L^2(X \rtimes \Gamma)$,
        \item It extends to an isometric $*$-isomorphism $A \rtimes_r \Gamma \rightarrow C_r^*(X \rtimes \Gamma)$.
    \end{enumerate}
\end{prop}
\begin{proof}
    Let $f,g \in C_c(\A_\beta)$ and $x \in X$, then 
    \begin{align*}
    \langle \Delta(f),\Delta(g) \rangle_{L^2(X \rtimes \Gamma)}(x) &= \sum\limits_{x \in D_t} \overline{\Delta(f)(x,t)}\Delta(g)(x,t)     
    = \sum\limits_{x \in D_t} \overline{f(t)(x)}g(t)(x).
    \end{align*}
    On the other hand, 
    \begin{align*}
    \langle f,g \rangle_{\ell^2(\A_\beta)}(x) = \left( \sum\limits_{t \in \Gamma} \overline{f(t)}g(t)  \right) (x)   
    = \sum\limits_{x \in D_t} \overline{f(t)(x)}g(t)(x).
    \end{align*}
    This proves (i). To prove (ii), note that $\Delta: \ell^2(\A_\beta) \rightarrow L^2(X \rtimes \Gamma)$ intertwines $\Lambda: C_c(\A_\beta) \rightarrow B(\ell^2(\A_\beta))$ with $\lambda: C_c(X \rtimes \Gamma) \rightarrow B(L^2(X \rtimes \Gamma)).$
\end{proof}
For $t \in \Gamma$, we write $\chi_t$ for the characteristic function on the set $D_t \times \{t\}$ and if $a \in C_0(D_t)$, then 
$$
\Delta(a\delta_t) = \chi_t * a .
$$
Thus we can rewrite Corollary \ref{emb Fell version} as the following.
\begin{prop}\label{emb Partial Actions}
There is an injective $*$-homomorphism $C_r^*(X \rtimes \Gamma) \rightarrow C_r^*(\Gamma) \otimes C_r^*(C \rtimes \Gamma)$ such that for all $t \in \Gamma$ and $a \in C_0(D_t)$, we have
    $$
\lambda(\chi_t * a) \mapsto \lambda^\Gamma_t \otimes \lambda(\chi_t * a).
    $$    
\end{prop}
    \subsection{Weak Amenability}
Let $\alpha = (\{D_t\}_{t \in \Gamma},\{\alpha_t\}_{t \in \Gamma}$) be a partial action of a discrete group $\Gamma$ on a locally compact Hausdorff space $X$ and $\G := X \rtimes \Gamma$ be the induced partial action groupoid and fix a quasi-invariant measure with full support $\mu$. Assume further that each $D_t$ is compact, thus $\G$ is an ample groupoid. In this section, we will prove that
$$
\cb(\G) = \cb(C_r^*(\G), C_0(X)).
$$
We already know $\cb(\G) \geq \cb(C_r^*(\G), C_0(X))$, so we only need to prove the other inequality. We write $A:= C_0(X)$ and given a completely bounded $A$-linear operator $T: C_r^*(\G) \rightarrow C_r^*(\G) $ we start by showing how to build a function $\varphi_T: \G \rightarrow \mathbb C$.
\begin{prop}\label{bebepequeno}
    Let $T:C_r^*(\G) \rightarrow C_r^*(\G)$ be an $A$-linear bounded operator and define for $(x,t) \in \G$
    \begin{equation}\label{function yeahhhhhh}
     \varphi_T(x,t) = E(\chi_t^* * T\chi_t)(x),   
    \end{equation}
    then:
    \begin{enumerate}[(i)]
    \itemsep0pt
        \item $\varphi_T \in C_b(\G)$ with $\| \varphi_T \|_\infty \leq \|T\|$;
        \item If $T$ has finite $A$-rank and takes values in $\Span_A\{\chi_t: t \in \Gamma\}$, then $\varphi_T \in C_c(\G)$ and $\|\varphi_T\|_{M_0A(\G)} \leq \|T\|_\text{cb}$.
    \end{enumerate}
\end{prop}
\begin{proof}
    For (i), we compute
    \begin{align*}
    |E( \chi_s^* *  T\chi_s)(x)| \leq \|\chi_s^* * T\chi_s \|_r  \leq \|T\|\|\chi_s\|_r^2 = \|T\|.
    \end{align*}
    It remains to prove that $\varphi_T$ is continuous. Since $ \G = \bigsqcup D_t \times \{t\} $ we prove that for all $t \in \Gamma$, $\varphi_T^t: D_t \rightarrow \mathbb C$, $\varphi_T^t(x) = E(\chi_t^* * T\chi_t)(x)$ is continuous. But in this case, $\varphi_T^t$ is just $E(\chi_t^* * T\chi_t)|_{D_t}$ which is continuous.
    For (ii), it suffices to show for the case where $T$ has $A$-rank one and takes values in $\Span_A\{\chi_{t_0}\}$ fixed $t_0$. This means there exists an $f \in C_r^*(\G)$ such that $T = \Theta_{\chi_{t_0},f}$. Let $(x,t) \in \G$, we compute
    $$
\varphi_T(x,t) = E(\chi_t^* * T\chi_t)(x) = E(\chi_t^* * \chi_{t_0} * \langle f,\chi_t \rangle_A )(x) = \chi_{t_0}(x,t)\langle f,\chi_t \rangle(x).
    $$
    Since $t_0$ is compact, we conclude $\varphi_T \in C_c(\G)$. We now prove the norm inequality. Let $V:L^2(\G,\nu^{-1}) \rightarrow \ell^2(\Gamma) \otimes L^2(\G,\nu^{-1})$ be the unitary operator $\chi_t * a \mapsto \delta_t \otimes \chi_t * a$ and $\tau: C_r^*(\G) \rightarrow C_r^*(\Gamma) \otimes C_r^*(\G)$ be the coaction of Proposition \ref{emb Partial Actions}, define $S: C_r^*(\G) \rightarrow C_r^*(\G)$ by 
    $$
S(f) = V^*(1 \otimes T)\tau(f)V,
    $$
    so if $f = \chi_t * a$ for some $s \in \Gamma$, $a \in C(D_t)$ we have
    \begin{equation}\label{operator S}
    S(\chi_s * a) = V^*(1 \otimes T)\tau(\chi_s * a)V = V^*(\lambda^\Gamma_s \otimes T(\chi_s * a))V  
    \end{equation}
    Our goal is to show that $S = M_{\varphi_T}$, for this, we keep the assumption that $T = \Theta_{\chi_{t_0},f}$. Using equation \eqref{operator S}, we get
    \begin{align*}
    S(\chi_s * a) = V^*(\lambda^\Gamma_s \otimes T(\chi_s * a))V = V^*(\lambda^\Gamma_s \otimes \lambda(\chi_{t_0} * \langle f,\chi_s \rangle_A ) )V     
    \end{align*}
    Let $s,t,u \in \Gamma$ and $a \in C(D_s),b \in C(D_t), c \in C(D_u)$ we compute
    \begin{align*}
        \langle S(\chi_s * a)(\chi_t * b), \chi_u * c \rangle &= \left\langle (1 \otimes T)\tau(\chi_{s} * a)(\delta_t \otimes \chi_t * b),\delta_u \otimes \chi_u * c \right\rangle \\
        &= \left\langle \lambda^\Gamma_s\delta_t \otimes \lambda(\chi_{t_0} * \langle f,\chi_s \rangle_A * a) (\chi_t * b),\delta_u \otimes \chi_u * c \right\rangle \\
        &= \langle \delta_{st},\delta_u \rangle \langle \chi_{t_0t} * d,\chi_u * c \rangle,
    \end{align*}
    where $d \in C(D_{t_0t})$ is such that $(\chi_{t_0} * \langle f,\chi_s \rangle_A )* (\chi_t * b) = \chi_{t_0t} * d$, that is,
    $$
    d = \beta_{t_0}\left(\beta_{t_0^{-1}}(a\langle f,\chi_s \rangle_A)b\right).
    $$
    On the other hand we compute 
    \begin{align*}
        \varphi_T(x,t)(\chi_s * a)(x,t) &= \begin{cases}
            \varphi_T(x,s)a(x) \quad &\text{ if } s=t, \\
            0 &\text{ otherwise.}
        \end{cases} \\
        &= \begin{cases}
            E(\chi_s^* * T\chi_s)(x) \quad &\text{ if } s=t, \\
            0 &\text{ otherwise.}
        \end{cases}\\
        &= \begin{cases}
            E(\chi_{s^{-1}t_0})(x) \langle f,\chi_s \rangle_A(x) a(x) \quad &\text{ if } s=t, \\
            0 &\text{ otherwise.}
        \end{cases}
    \end{align*}
    $E(\chi_{s^{-1}t_0})$ is the restriction to $X$ of $\chi_{s^{-1}t_0}$, which is only non-zero if $s = t_0$. We conclude
    \begin{equation}
        M_{\varphi_T}(\chi_s * a) = \begin{cases}
            \chi_{t_0}* a\langle f,\chi_{t_0} \rangle_A   \quad &\text{ if } s=t_0, \\
            0 &\text{ otherwise.}
        \end{cases}
    \end{equation}
    Using the equation above we are able to calculate
    \begin{align*}
        \langle M_{\varphi_T}(\chi_s * a)(\chi_t * b),(\chi_u * c) \rangle &= \langle \delta_{t_0},\delta_s \rangle \langle \lambda(\chi_{t_0}* a\langle f,\chi_{t_0} \rangle_A)(\chi_t*b),\chi_u * c \rangle \\
        &=  \langle \delta_{t_0},\delta_s \rangle \langle \chi_{t_0t} * d,\chi_u * c \rangle
    \end{align*}
    The left side is only non-zero if $t_0t = u$, in this case we have $t_0 = s \iff st = u$ and we conclude 
    $$
    \langle M_{\varphi_T}(\chi_s * a)(\chi_t * b),(\chi_u * c) \rangle = \langle \delta_{st},\delta_u \rangle \langle \chi_{t_0t} * d,\chi_u * c \rangle.
    $$
    This proves that $S = M_{\varphi_T}$. Moreover, the general case where $T$ has finite $A$-rank follows from the linearity of the inner product. We finish the proof by noting that due to the definition of $S$,
    $$
\|M_{\varphi_T}\|_{\text{cb}} = \|S\|_{\text{cb}} \leq \|T\|_{\text{cb}}.
    $$
\end{proof}
With this out of the way, we are almost there. Suppose now that $T_i: C_r^*(\G) \rightarrow C_r^*(\G)$ is a net testifying the CBAP of the pair ($C_r^*(\G),C_0(X))$, of course, it is not necessary that each $T_i$ satisfies condition (ii) of Proposition \ref{bebepequeno}, that is, $T_i$ might not take values in $\Span_A\{\chi_t: t \in \Gamma\}$. We now present a proposition that will allow us to build a new net, for which condition (ii) is met. It is a quasi Cartan pair version of \cite[Lemma 11.4]{buss2024fourier}.
\begin{lema}\label{rotatelemma}
   Let $(B,A)$ be a quasi Cartan pair, $T:B \rightarrow B$ be a finite $A$-rank map and $B_0 \subseteq B$ be a dense subspace. Then, there exists a sequence $(T_n)$ of finite $A$-rank maps $T_n : B \rightarrow B_0 \subseteq B$ such that $\|T - T_n\|_\text{cb} < \frac{1}{n}$. 
\end{lema}
\begin{proof}
Since $T$ has finite $A$-rank there are $g_i,h_i \in B, \, i = 1,...,N$ such that  
$$
T = \sum\limits_{i=1}^N \Theta_{g_i,h_i}
$$
Using Proposition \ref{decomposablerankone} we compute
$$
\|T\|_{\text{cb}} =  \sum\limits_{i=1}^N \| \Theta_{g_i,h_i} \|_\text{cb} \leq \sum\limits_{i=1}^N \|g_i\|\|h_i\| 
$$
For each $i=1,...,N$ and $n \in \nat$, use the density of $B_0$ to choose $g_{i,n} \in B_0$ such that $\left\|g_i - g_{i,n} \right\| < 1/(nN\max_i \|h_i\|)$  and let 
$$
T_n := \sum\limits_{i=1}^N \Theta_{g_i,h_{i,n}} = \sum\limits_{i=1}^N g_{i,n} * \langle h_{i}, \cdot \rangle_A \in B.
$$
Since $T_n$ has also finite $A$-rank, by the above argument
$$
\|T - T_n\|_{\text{cb}} \leq \sum\limits_{i=1}^N \|g_i\| \, \|h_i - h_{i,n}\| < \frac{1}{n}.
$$
\end{proof}
We now have the two main ingredients.
\begin{theorem}\label{theorempartialactions}
    Let $\alpha = (X,\Gamma,\{D_t\}_{t \in \Gamma},\{\alpha_t\}_{t \in \Gamma}$) be a partial action of a discrete group $\Gamma$ on a locally compact Hausdorff space $X$. Assume further that each $D_t$ is compact. Let $\G = X \rtimes \Gamma$ be the ample measured groupoid arising from this partial action. Then
    $$
\cb(\G) = \cb(C_r^*(\G),C_0(X)).
    $$
\end{theorem}
\begin{proof}
    By Proposition \ref{weakinequality}, we know that $\cb(\G) \geq \cb(C_r^*(\G),C_0(X)).$ For the opposite inequality, let $(T_i)_{i \in I}$ be a net testifying the CBAP of the quasi Cartan pair $(C_r^*(\G),C_0(X))$. Applying Lemma \ref{rotatelemma} to each $i$, we obtain a new net of finite $A$-rank completely bounded maps $T'_{i,n}$ such that $T'_{i,n} \xrightarrow{\text{SOT}} \id_{C_r^*(\G)}$ (using the product directed set) with the added benefit that each $T'_{i,n}$ takes values in $\Span_A\{\chi_t: t \in \Gamma\}$, however it does not satisfy the norm requirements. This is easily fixed by letting 
    $$
T_{i,n} := \frac{C}{C+1/n}T'_{i,n}.
    $$
    For every $i \in I$ and $n \in \nat$,
    \begin{align*}
     \|T_{i,n}\|_{\text{cb}} =  \frac{C}{C+1/n} \|T'_{i,n}\|_{\text{cb}} <   \frac{C}{C+1/n} \left( C + \frac{1}{n} \right) = C.
    \end{align*}
    We conclude that the net $(T_{i,n})$ testifies the conditions of the CBAP of $(C_r^*(\G),C_0(X))$. Using Proposition \ref{bebepequeno}, for every $i \in I$ and $m \in \nat$ we obtain $\varphi_{T_{i,n}} \in C_c(\G) $ such that $\|\varphi_{T_{i,n}}\|_{M_0A(\G)} \leq C$. Let $K \subseteq \G$ be compact and cover it with 
    $$
    \{D_{t_j} \times \{t_j\}: t_j \in \Gamma \text{ for } j = 1,...,N \},
    $$ 
    Let $\varepsilon > 0$ and recalling $T_{i,n} \xrightarrow{\text{SOT}} \id_{C_r^*(\G)}$ choose for each $j=1,...,N$ an $M_j \in I \times \nat$ such that $(i,n) > M_j$ implies $\|T_{i,n}\chi_{t_j} - \chi_{t_j} \|_r < \varepsilon$. If we let $M := \max_{j=1,...,N} M_j$, $(i,n) > M$, $(x,t_k) \in K$ we have
    \begin{align*}
      |1-\varphi_{i,n}(x,t_k)| &= | E(\chi_{t_k}^* * \chi_{t_k})(x) - E(\chi_{t_k}^* * T_{i,n}\chi_{t_k})(x) | \\
      &= |E( \chi_{t_k}^* * ( \chi_{t_k} - T_{i,n}\chi_{t_k}) )(x) | \\
      &\leq \|\chi_{t_k}^* * (\chi_{t_k} - T_{i,n}\chi_{t_k})\|_r \\
      &\leq \|\chi_{t_k} - T_{i,n}\chi_{t_k}\|_r \\
      &< \varepsilon,
    \end{align*}
    proving that $\varphi \rightarrow 1$ uniformly on compacts. We conclude that $\varphi_{i,n}$ testifies the weak amenability of $\G$ and that
    $$
    \cb(\G) = \cb(C_r^*(\G),C_0(X)).
    $$
\end{proof}
\section{Weak Amenability of Discrete Groupoids}
Let $\G$ be a discrete groupoid and suppose $\mu$ is a quasi-invariant measure on $X := \gob$. Since $\G$ is discrete, 
$$
\supp(\mu) := \{ x \in X: \mu(\{x\}) > 0 \}
$$
so if $\mu'$ is another quasi-invariant measure on $\gob$ and $\supp(\mu') = \supp(\mu)$ then $\mu'$ and $\mu$ are equivalent. Thus if we are given a discrete measured groupoid $(\G,\mu)$ by passing to $(\G|_\Omega,\mu)$, where $\Omega := \supp(\mu)$ we can assume $\mu$ has full support. Our goal will be to show that
$$
\cb(\G) = \cb(C_r^*(\G),C_0(X)).
$$
We do this by introducing Fell's absorption principle for discrete groupoids and with it proving there exists a normal $*$-embedding $\pi: \vN(\G) \rightarrow \vN(\G) \, \overline{\otimes} \, \vN(\G)$, this is done in Section 8.1. Note that we write simply $\vN(\G)$ instead of $\vN(\G,\mu)$ due to our previous discussion, this algebra only depends on the support of $\mu$, which we assume to be $X$. In Section 8.2, we prove the main result using the methods in Haagerup's paper \cite{haagerup2016group}.
    \subsection{Fell Absorption for Discrete Groupoids}

Unlike what happens for groups, Fell absorption for discrete groupoids is given by partial isometries instead of unitary operators.
\begin{defi}
    Let $H,K$ be Hilbert spaces, a bounded linear operator $W: H \rightarrow K $ is said to be a \textit{partial isometry} if $W = WW^*W$. We call $\Init(W) := (\ker W)^\perp$ its \textit{initial subspace} and $\Fin(W) := \im W$ its \textit{final subspace}.
\end{defi}
\begin{example}
    Let $\gamma \in \G$, we prove that $\lambda_\gamma: \L \rightarrow \L$ is a partial isometry. Let $f \in \L$, then
    \begin{align*}
        \lambda_\gamma\lambda_\gamma^*\lambda_\gamma(f)(\beta) = \begin{cases}
        \lambda_\gamma^*\lambda_\gamma(f)(\gamma^{-1}\beta), \quad &\text{if } \beta \in \G^{r(\gamma)} \\
        0 &\text{otherwise}
        \end{cases} = \begin{cases}
        \lambda_\gamma(f)(\gamma\gamma^{-1}\beta), \quad &\text{if } \beta \in \G^{r(\gamma)} \\
        0 &\text{otherwise}
        \end{cases} 
        = \lambda_\gamma(f)(\beta).
    \end{align*}
    Moreover, it is easy to see that the initial subspace of $\lambda_\gamma$ is $L^2\left(\G^{d(\gamma)},\nu^{-1}\right)$ and the final subspace $L^2(\G^{r(\gamma)},\nu^{-1})$.
\end{example}
\begin{example}
    More generally, it is easy to see that if $\pi:C_c(\G) \rightarrow B(H)$ is a $*$-representation, then $\pi(\gamma)$ is a partial isometry.
\end{example}
\begin{prop}(Fell absorption for discrete groupoids)\label{discrete Fell}
Let $(\G,\mu)$ be a discrete measured groupoid, where $\mu$ has full support and $\pi: C_c(\G) \rightarrow B(H)$ a $*$-representation. Then there exists a partial isometry $W_\lambda: H \otimes \L \rightarrow H \otimes \L$ such that for every $\gamma \in \G$,
$$
W_\pi(1 \otimes \lambda_\gamma)W^*_\pi = \pi(\gamma) \otimes \lambda_\gamma.
$$
\end{prop}
\begin{proof}
    We define $W_\pi: H \otimes \L \rightarrow H \otimes \L$ by $\xi \otimes \delta_\beta \mapsto \pi(\beta)\xi \otimes \delta_\beta$, for all $\xi \in H$ and $\beta \in \G$. We readily see its a partial isometry
    \begin{align*}
     W_\pi W^*_\pi W_\pi(\xi \otimes \delta_\beta) = W_\pi W_\pi^*(\pi(\beta)\xi \otimes \delta_\beta) = W_\pi(\pi(\beta^{-1})\pi(\beta)\xi \otimes \delta_\beta) = W_\pi(\xi \otimes \delta_\beta).   
    \end{align*}
    Now for $\gamma \in \G$, we compute 
    \begin{align*}
    W_\pi(1 \otimes \lambda_\gamma)W_\pi^*(\xi \otimes \delta_\beta) = W_\pi(1 \otimes \lambda_\gamma) (\pi(\beta^{-1})\xi \otimes \delta_\beta) &=  \begin{cases}
        W_\pi(\pi(\beta^{-1})\xi \otimes \lambda_{\gamma\beta}), \quad &\text{ if } r(\beta) = d(\gamma) \\
        0 &\text{ otherwise}
        \end{cases} \\
       &= \begin{cases}
        \pi(\gamma\beta\beta^{-1})\xi \otimes \lambda_{\gamma\beta}, \quad &\text{ if } r(\beta) = d(\gamma) \\
        0 &\text{ otherwise}
        \end{cases} \\
        &= \begin{cases}
        \pi(\gamma)\xi \otimes \lambda_{\gamma\beta}, \quad &\text{ if } r(\beta) = d(\gamma) \\
        0 &\text{ otherwise} 
        \end{cases} \\
        &= \pi(\gamma) \otimes \lambda_\gamma(\xi \otimes \delta_\beta).
    \end{align*}
\end{proof}
We will want to create an embedding $\pi: \vN(\G) \rightarrow \vN(\G,\mu) \overline{\otimes} \vN(\G)$, but before, a couple of lemmas.
\begin{lema}\label{embtensor}
    Let $H,K$ be Hilbert spaces and $\{e_i\}$ an orthonormal basis for $K$. Let $W:H \otimes K \rightarrow H \otimes K$ be a partial isometry with initial subspace $M$. Assume that for all $i \in I$, there exist a non-zero $v \in H$ such that $v \otimes e_i \in M$. Let $A \subseteq B(H)$ be a $*$-subalgebra such that for all $T \in A$, $(1 \otimes T)(M) \subseteq M$, then the $*$-homomorphism
    $$
  \pi:  A \ni T \mapsto W(1 \otimes T)W^* \in B(H)
    $$
    is injective.
\end{lema}
\begin{proof}
    Assume $W(1 \otimes T)W^*=0$, let $i \in I$ and $v$ such that $v \otimes e_i \in M$ so that $W(v \otimes e_i) \in N$, then 
    $$
W(1 \otimes T)W^*(W(v \otimes e_i)) = W(1 \otimes T)(v \otimes e_i) = W(v \otimes Te_i) = 0.
    $$
    Since $1 \otimes T$ is invariant for $M$ in which $W$ is isometric we conclude that 
    $v \otimes Te_i = 0$, recalling that $v \neq 0$ this implies $Te_i = 0$, since $i$ is arbitrary, $T=0$ as desired.
\end{proof}
\begin{lema}\label{invariance}
Let $\lambda: C_c(\G) \rightarrow B(\L)$ be the left regular representation of $\G$ and let $W_\lambda:\L \otimes \L \rightarrow \L \otimes \L$ the partial isometry of Proposition \ref{discrete Fell}, then the initial subspace of $W_\lambda$ is 
    $$
\text{Init}(W_\lambda) = \cls \{ \delta_\alpha \otimes \delta_\beta : \alpha \in \G^{d(\beta)} \text{ and } \beta \in \G \}.
    $$
Moreover, for every $\gamma \in \G$,
$
1 \otimes \lambda_\gamma(\Init(W_\lambda)) \subseteq \Init (W_\lambda).
$
\end{lema}
\begin{proof}
It is well known that $\Init(W_\lambda) = \Fin(W_\lambda^*) = \im(W_\lambda^*)$, so let $\alpha,\beta \in \G$ and compute 
\begin{align*}
    W_\lambda^*(\delta_\alpha \otimes \delta_\beta) = \lambda_{\beta^{-1}}\delta_\alpha \otimes \delta_\beta = \begin{cases}
        \delta_{\beta^{-1}\alpha} \otimes \delta_\beta, \quad &\text{ if } \alpha \in \G^{r(\beta)}, \\
        0  &\text{ otherwise.}
        \end{cases}
\end{align*}
In all cases, we have that $\lambda_{\beta^{-1}}\delta_\alpha \otimes \delta_\beta \in \{ \delta_\alpha \otimes \delta_\beta : \alpha \in \G^{d(\beta)} \text{ and } \beta \in \G \}$, proving the first part. Now, let $\delta_\alpha \otimes \delta_\beta \in \Init(W_\lambda) $, then 
$$
1 \otimes \lambda_\gamma(\delta_\alpha \otimes \delta_\beta) = \begin{cases}
         \delta_\alpha \otimes \delta_{\gamma\beta}, \quad &\text{ if } \beta \in \G^{d(\gamma)}, \\
        0  &\text{ otherwise.}
        \end{cases}
$$
By definition $\alpha \in \G^{d(\beta)} = \G^{d(\gamma\beta)}$ so $\delta_\alpha \otimes \delta_{\gamma\beta} \in \Init(W_\lambda)$, as desired.
\end{proof}
We can now prove our desired embedding.
\begin{theorem}\label{discrete embedding}
    There exists a normal $*$-embedding $\pi: \vN(\G) \rightarrow \vN(\G) \overline{\otimes} \vN(\G) $ such that for all $\gamma \in \G$,
    $$
\pi(\lambda_\gamma) = \lambda_\gamma \otimes \lambda_\gamma.
    $$
\end{theorem}
\begin{proof}
    For $f \in C_c(\G)$, we define 
    $$
\pi(f) = W_\lambda(1 \otimes \lambda(f))W^*_\lambda,
    $$
    it is just the map of Lemma \ref{embtensor}, so to show it is injective, we just need to prove the conditions of said lemma. We have $H = K = \L$, $W = W_\pi$ and that $\{\delta_\gamma: \gamma \in \G\}$ is an orthonormal basis for $\L$. The initial subspace of $W_\lambda$ has been computed in Lemma \ref{invariance} to be 
    $$
\text{Init}(W_\lambda) = \cls \{ \delta_\alpha \otimes \delta_\beta : \alpha \in \G^{d(\beta)} \text{ and } \beta \in \G \}
    $$
    and note that for every $\gamma \in \G$, $\delta_\gamma \otimes \delta_\gamma \in \Init(W_\lambda)$. In the conditions of Lemma \ref{embtensor}, we have $A = \lambda(C_c(\G)) \subseteq B(H)$ and the invariance of $\Init(W)$ for $(1 \otimes \lambda(f))$ follows from the second part of Lemma \ref{invariance}. Moreover, due to Proposition \ref{discrete Fell}, $\pi(\lambda_\gamma) = \lambda_\gamma \otimes \lambda_\gamma$. To finish the proof, we note that $\pi$ extends to $\pi:C_r^*(\G) \rightarrow B(\L \otimes \L)$ and by definition of the minimal tensor product, it actually takes values in $C_r^*(\G) \otimes C_r^*(\G)$. Since $\pi$ is spatially implemented (by $W_\lambda$) then by \cite[Theorem 5.2]{takesaki2001theory} it extends to a normal $*$-embedding $\pi: \vN(\G) \rightarrow \vN(\G) \overline{\otimes} \vN(\G)$.
\end{proof}
    \subsection{Weak Amenability}
Recall the map $\pi: \vN(\G) \rightarrow \vN(\G) \overline{\otimes} \vN(\G)$ of Theorem \ref{discrete embedding}. The von Neumann algebra $\vN(\G)$ admits a normal faithful weight $\tau$ (see \cite[Lemma 2.5]{anantharaman2013haagerup}) given by 
$$
\tau(T) = \langle \chi_X, T\chi_X \rangle_{\L}, \quad T \in \vN(\G).
$$
The one-parameter modular automorphism group relative to $\tau$ and $\chi_X$ has been computed in \cite{hahn1978regular} and we present it as in \cite[p.6]{anantharaman2013haagerup}: for $T \in \vN(\G)$, we have
$$
\sigma_t(T) = D^{it}TD^{-it},
$$
where $D^{it}$ acts on $\L$ by pointwise multiplication and defines a unitary operator.
\begin{prop}
    For all $t \in \mathbb R$, $\sigma_t|_{\im(\pi)} = \id_{\vN(\G)}$. 
\end{prop}
\begin{proof}
Let $\gamma \in \G$, then we compute for $\beta \in \G$
\begin{align*}
\sigma_t|_{\im(\pi)}(\delta_\beta) =  D^{it}\lambda_\gamma D^{-it}(\delta_\beta) =  D^{it}\lambda_\gamma (D^{-it}(\beta)\delta_\beta) &=  \begin{cases}
        D^{it}(D^{-it}(\beta)\delta_{\gamma\beta}), \quad &\text{ if } r(\beta) = d(\gamma) \\
        0 &\text{ otherwise}
        \end{cases} \\
        &= \begin{cases}
        \delta_{\gamma\beta}, \quad &\text{ if } r(\beta) = d(\gamma) \\
        0 &\text{ otherwise}
        \end{cases} \\
        &= \lambda_\gamma(\delta_\beta).
\end{align*}
Since $\{\lambda_\gamma:\gamma \in \G \}$ generates $\vN(\G)$, we are finished
\end{proof}
Now, by \cite[Proposition 4.3]{takesaki2003theory} the weight $\tau \otimes \tau$ is also semi-finite and faithful and we have 
$$
\sigma_t^{\tau \otimes \tau} = \sigma_t \otimes \sigma_t = \id_{\vN(\G)} \otimes \id_{\vN(\G)}, \quad \text{for all } t \in \mathbb R.
$$
This means $\sigma_t \otimes \sigma_t(\im\pi) = \im\pi$ so we can apply \cite[Theorem 4.2]{takesaki2003theory} to obtain a conditional expectation $\mathcal{E}: \vN(\G) \overline{\otimes} \vN(\G) \rightarrow \im \pi$ such that $\tau \otimes \tau \circ \mathcal{E} = \tau \otimes \tau$. We abstract our setting to prove the following lemma.
\begin{lema}\label{lemma weight}
    In the language of \cite{takesaki2003theory}, let $\varphi$ be a faithful semi-finite weight on a von Neumann algebra $\mathcal M$ and $\mathcal N$ be a von Neumann subalgebra of $\mathcal M$ such that the restriction $\varphi|_{\mathcal N}$ of $\varphi$ to $\mathcal N$ is semi-finite. Suppose there exists a conditional expectation $\mathcal{E}: \mathcal M \rightarrow \mathcal N$ verifying $\varphi = \varphi \circ \mathcal E$. Then, for each $x \in  \mathcal{N}$, $\mathcal E(x)$ is the unique element in $ \mathcal{N} $ such that for all $y,z \in \mathfrak n_\varphi \cap \mathcal{N}$,
    \begin{equation}\label{equalitywieght}
    \varphi(z^*\mathcal{E}(x)y) = \varphi(z^*xy).    
    \end{equation}
\end{lema}
\begin{proof}
Let $x \in \mathcal{N}$, we first prove uniqueness. We recall the setup as in the book \cite{takesaki2003theory}, we define the set $\mathfrak p_\varphi := \{ x \in \mathcal M_+ : \varphi(x) < \infty \} $ and $\mathfrak n_\varphi := \{ x \in \mathcal M : x^*x \in \mathfrak p_\varphi \}$. We also let $\mathfrak p_\varphi ' = \mathfrak p_{\varphi|_{\mathcal{N}}}$ and similarly $\mathfrak n_\varphi ' = \mathfrak n_{\varphi|_{\mathcal{N}}}$. It is easy to see that $\mathfrak p'_\varphi = \mathfrak p_\varphi \cap \mathcal{N}$ and $\mathfrak n'_\varphi = \mathfrak n_\varphi \cap \mathcal{N}$. Since $\varphi|_{\mathcal{N}}$ is semi-finite $\mathfrak n'_\varphi$ is non-empty and we define an inner product on it by (recall that $\varphi$ is faithful) 
$$
\langle y,z \rangle_\varphi := \varphi(z^*y), \quad \text{ for all } y,z \in \mathfrak n_\varphi'.
$$
We denote by $H_\varphi$ the Hilbert space obtained via completion. $\mathcal{N}$ acts on $\mathfrak n_\varphi'$ by multiplication and we can extend this to $H_\varphi$, obtaining a representation $\pi_\varphi : \mathcal{N} \rightarrow B(H_\varphi)$, it is faithful since $\varphi$ is \cite[Proposition 1.4]{takesaki2003theory} and from now on we view $\mathcal{N} \subseteq B(H_\varphi)$. To prove that there exists a unique $T_x \in B(H_\varphi)$ satisfying $\varphi(z^*T_xy) = \varphi(xy)$, we use Riesz representation theorem, indeed, $S_x:H_\varphi \times H_\varphi \rightarrow \mathbb{C}$ defined as $ S_x(y,z) = \varphi(z^*xy)$ is a bounded sesquilinear functional meaning there exists a unique $T_x \in B(H_\varphi)$ such that
$$
\langle T_xy,z \rangle_\varphi = S_x(y,z) \iff \varphi(z^*T_xy) = \varphi(z^*xy), \quad \text{ for all } y \in H_\varphi, 
$$
proving uniqueness. Then, the fact that $T_x = \mathcal{E}(x) \in \mathcal N \subseteq B(H_\varphi)$ follows from the equation $\varphi = \varphi \circ \mathcal{E}$, indeed if $y,z \in \mathfrak n'_\varphi $ we have 
$$
\varphi(z^*xy) = \varphi(\mathcal{E}(z^*xy)) = \varphi(z^*\mathcal{E}(x)y).
$$
\end{proof}
\begin{prop}\label{condgroupobundles}
    The conditional expectation $\mathcal{E}:\vN(\G) \overline{\otimes} \vN(\G) \rightarrow \im \pi$ satisfies for all $\alpha,\beta \in \G$ 
    $$
\mathcal{E}(\lambda_\alpha \otimes \lambda_\beta) = \begin{cases}
    \lambda_\alpha \otimes \lambda_\alpha \quad &\text{ if } \alpha = \beta, \\
    0 &\text{ otherwise.}
\end{cases}
    $$
\end{prop}
\begin{proof}
We make use of Lemma \ref{lemma weight}, which falls in our setting with $\mathcal{M} = \vN(\G) \overline{\otimes} \vN(\G)$, $\mathcal{N} = \im(\pi)$ and $\varphi = \tau \otimes \tau$. Let $\gamma_1,\gamma_2 \in \G$, then it is easy to see that $\lambda_{\gamma_i} \otimes \lambda_{\gamma_i} \in \mathfrak n_\varphi \cap \mathcal N$, we compute for $\alpha,\beta \in \G$ 
\begin{align*}
(\lambda_{\gamma_2} \otimes \lambda_{\gamma_2})^*(\lambda_{\alpha} \otimes \lambda_{\beta})(\lambda_{\gamma_1} \otimes \lambda_{\gamma_1}) &= \begin{cases}
        (\lambda_{\gamma_2^{-1}} \otimes \lambda_{\gamma_2^{-1}})(\lambda_{\alpha\gamma_1} \otimes \lambda_{\beta\gamma_1}) \quad &\text{ if } r(\gamma_1) = d(\alpha) \\
        &\text{ and } r(\gamma_1) = d(\beta), \\
        0 &\text{ otherwise.} 
        \end{cases} \\
        &= \begin{cases}
        \lambda_{\gamma_2^{-1}\alpha\gamma_1} \otimes \lambda_{\gamma_2^{-1}\beta\gamma_1} \quad &\text{ if } r(\gamma_1) = d(\alpha), r(\gamma_1) = d(\beta) \\
        &\text{ and }r(\gamma_2) = r(\alpha), r(\gamma_2) = r(\beta), \\
        0 &\text{ otherwise.}
        \end{cases}    
\end{align*}
For the case where the above expression is non-zero we compute
\begin{align*}
\tau \otimes \tau \left( \lambda_{\gamma_2^{-1}\alpha\gamma_1} \otimes \lambda_{\gamma_2^{-1}\beta\gamma_1} \right) = \left\langle \chi_X,\delta_{\gamma_2^{-1}\alpha\gamma_1} \right\rangle  \left\langle \chi_X,\delta_{\gamma_2^{-1}\beta\gamma_1} \right\rangle,   
\end{align*}
The right hand side is non-zero if $x_\alpha := \gamma_2^{-1}\alpha\gamma_1 \in X$ and $x_\beta := \gamma_2^{-1}\beta\gamma_1$, note that we necessarily have that $x_\alpha = x_\beta =: x \in X$, this implies 
$$
\gamma_2^{-1}\alpha\gamma_1 = \gamma_2^{-1}\beta\gamma_1 \iff \alpha = \beta.
$$
Recalling that by definition $\mathcal{E}|_{\im(\pi)} = \id_{\im(\pi)}$ we conclude 
$$
\mathcal{E}(\lambda_\alpha \otimes \lambda_\beta) = \begin{cases}
    \lambda_\alpha \otimes \lambda_\alpha \quad &\text{ if } \alpha = \beta, \\
    0 &\text{ otherwise.}
\end{cases}
    $$
\end{proof}
Before we move on to our main theorem, we recall that if $T: C_r^*(\G) \rightarrow C_r^*(\G)$ is a finite $C_0(X)$-rank operator, then we can "rotate" it via Lemma \ref{rotatelemma} to obtain a sequence $T_n: C_r^*(\G) \rightarrow C_r^*(\G)$ of finite $C_0(X)$-rank operators takes values in a dense subset $B_0 \subseteq C_r^*(\G)$. Since $\G$ is discrete, we take 
$$
B_0 := \Span\{\lambda_\gamma: \gamma \in \G\}.
$$
We now show how to obtain a function $\varphi_T$ from an $C_0(X)$-linear bounded operator $T:C_r^*(\G) \rightarrow C_r^*(\G)$ similar to what we did in Section 7.
\begin{prop}\label{bebepequenodiscrete}
    Let $T:C_r^*(\G) \rightarrow C_r^*(\G)$ be a $C_0(X)$-linear bounded operator and define for $\gamma \in \G$
    \begin{equation}\label{function yeahhhhhh}
     \varphi_T(\gamma) = E(\lambda_{\gamma^{-1}} * T\lambda_\gamma)(d(\gamma)),   
    \end{equation}
    then:
    \begin{enumerate}[(i)]
    \itemsep0pt
        \item $\varphi_T \in C_b(\G)$ with $\| \varphi_T \|_\infty \leq \|T\|$;
        \item If $T$ has finite $C_0(X)$-rank and takes values in $\Span\{\lambda_\gamma: \gamma \in \G \}$, then $\varphi_T \in C_c(\G)$. Moreover, $\|\varphi_T\|_{M_0A(\G)} \leq \|T\|_\text{cb}$.
    \end{enumerate}
\end{prop}
\begin{proof}
(i) is proven very similarly to Proposition \ref{bebepequeno}. For (ii), it suffices to assume $T$ is a $C_0(X)$-rank one operator taking values in $\Span\{\lambda_{\gamma}\}$ for a $\gamma \in \G$, so there exists an $f \in C_r^*(\G)$ such that $T = \Theta_{\lambda_{\gamma},f}$. Let $\beta \in \G$, then there exists an $a_\beta \in \mathbb C$ such that $T\lambda_\gamma =a_\beta\lambda_{\gamma} $ and we have
\begin{align*}
\varphi_T(\beta) = E(\lambda_{\beta^{-1}} * T\lambda_\beta)(d(\beta)) = E(\lambda_{\beta^{-1}} * a_\beta\lambda_{\gamma})(d(\beta))    
\end{align*}
which is only non-zero if $\beta = \gamma$, so $\varphi_T \in C_c(\G)$. To prove the norm inequality, define an operator $S: C_r^*(\G) \rightarrow C_r^*(\G)$ by 
$$
S = \pi^{-1} \circ \mathcal{E} \circ( 1 \otimes T ) \circ \pi
$$
and for $\beta \in \G$, compute 
\begin{align*}
\pi^{-1} \circ \mathcal{E} \circ( 1 \otimes T ) \circ \pi(\lambda_\beta) =  \pi^{-1} \circ \mathcal{E}( \lambda_\beta \otimes T \lambda_\beta )  =  \pi^{-1} \circ \mathcal{E} ( \lambda_\beta \otimes a_\beta\lambda_\gamma ) &= \pi^{-1}(a_\beta\mathcal{E}( \lambda_\beta \otimes \lambda_\gamma )).
\end{align*}
By Proposition \ref{condgroupobundles}, we conclude that 
$$
S(\lambda_\beta) = \begin{cases}
    a_\beta\pi^{-1}(\lambda_\gamma \otimes \lambda_\gamma) \quad &\text{ if } \gamma = \beta, \\
    0 &\text{ otherwise.}
\end{cases} = \begin{cases}
    a_\beta\lambda_\gamma  \quad &\text{ if } \gamma = \beta, \\
    0 &\text{ otherwise.}
\end{cases}
$$
On the other hand, 
\begin{align*}
    m_{\varphi}(\lambda_\beta) = \varphi(\beta)\lambda_\beta = E(\lambda_\beta^* * T\lambda_\beta)(d(\beta))\lambda_\beta &= E(\lambda_\beta^* * a_\beta\lambda_\gamma)(d(\beta))\lambda_\beta \\
    &=  \begin{cases}
    a_\beta\lambda_\gamma  \quad &\text{ if } \gamma = \beta, \\
    0 &\text{ otherwise.}
\end{cases} \\
   &= S(\lambda_\beta).
\end{align*}
We obtain $S = m_\varphi$ and since $\|S\|_{\text{cb}} \leq \|T_{n,m}\|_{\text{cb}} \leq C$ we conclude that $\|\varphi_{n,m}\|_{M_0A(\G)} \leq C$, finishing the proof.
\end{proof}
\begin{theorem}\label{cb discrete equality}
    Let $\G$ be a discrete groupoid, then 
    $$
\cb(\G) = \cb(C_r^*(\G),C_0(X)).
    $$
\end{theorem}
\begin{proof}
By Proposition \ref{weakinequality}, we know that $\cb(C_r^*(\G),C_0(X)) \leq  \cb(\G) $. Let $T_i$ be a net for the CBAP of the quasi Cartan pair $(C_r^*(\G),C_0(X))$ and $C>0$ such that $\sup_i\|T_i\|_{\text{cb}} \leq C$. Similarly to what we did in the proof of Theorem \ref{theorempartialactions}, we can assume without loss of generality that each $T_i$ takes values in $\Span\{\delta_\gamma : \gamma \in \G\}$. By Proposition \ref{bebepequenodiscrete}, we get functions $\varphi_{i} := \varphi_{T_i} \in C_c(\G)$ such that $\|\varphi_{i}\|_{M_0A(\G)} \leq C$. Since $\G$ is discrete, it simply remains to show that $\varphi_i \rightarrow 1$ pointwise. For each $\gamma \in \G$, choose $h_\gamma^{(i)} \in C_r^*(\G)$, where only a finite number of them are non-zero such that 
$$
T_i = \sum\limits_{\gamma \in \G} \Theta_{\lambda{\gamma},h_\gamma^{(i)}}
$$
For $\beta \in \G$, let $a_\gamma^{(i)} \in \mathbb C$ be such that $\Theta_{\lambda{\gamma},h_\gamma^{(i)}}(\lambda_\beta) = a_\gamma^{(i)}\lambda_{\gamma} $, then 
$$
\varphi_i(\beta) = E\left( \lambda_{\beta^{-1}} * \sum\limits_{\gamma \in \G} a_\gamma^{(i)}\lambda_{\gamma} \right) (d(\gamma)) =  a_\beta^{(i)}. 
$$
Since $T_i \xrightarrow{\text{SOT}} \id_{C_r^*(\G)}$, it can be proven that $a_\beta^{(i)} \rightarrow 1$, finishing the proof.
\end{proof}
\section{Open Questions and Further Work}
We finish this work by laying out some questions that arose in the process and what other research is possible in this theory. 
\begin{question}
Does the equality $\cb(\G) = \cb(C_r^*(\G),C_0(\gob))$ hold for all Hausdorff étale groupoids $\G$?    
\end{question}
This was the main goal of the paper, which unfortunately was not met. The extensive range of examples given in Sections 7 and 8 strongly suggest that the claim is true. The biggest impediment into proving this was the fact that there does not exist a strong enough Fell absorption principle that yields a "diagonal" coaction $\pi:C_r^*(\G) \rightarrow \A \otimes C_r^*(\G)$ for a suitable $C^*$-algebra $\A$. Nonetheless, a general Fell absorption principle does exist, done in \cite{buss2023approximation}.

\begin{question}
Does measurewise weak amenability imply topological amenability for all Hausdorff étale groupoids $\G$?       
\end{question}

Most likely yes. An argument similar to \cite[Theorem 3.3.7]{anantharaman2001amenable} where it is proven that topological  and measurewise amenability are equivalent might work. However, going deep into the analysis part of weak amenability was outside of the scope of this work.

\begin{question}
Does Cowling's Conjecture hold for Hausdorff étale groupoids?   
\end{question}
If $G$ is a locally compact group, Cowling's conjecture reads:
\begin{conjecture}
If $G$ is weakly amenable with $\cb(G) = 1$, then it has the Haagerup property.    
\end{conjecture}
The Haagerup property for Hausdorff étale groupoids $\G$ has been defined in \cite[Proposition 5.4]{kwasniewski2022haagerup}. Thus, one could pose the same question but substituting the locally compact group $G$ with the groupoid $\G$. Although, it will of course be harder to prove if the conjecture is true, if instead one wants to attempt to find a counter-example, the groupoid world might be a good start, similarly to what happened regarding the Baum-Connes conjecture.

\begin{question}
What about non Hausdorff étale groupoids?
\end{question}
The Fourier-Stieltjes and Fourier algebras can be defined just fine. However, a recent paper of Buss and Martínez \cite{buss2025essential} suggest that an \textit{essential} theory is more correct. To achieve this, we would have to adapt the definitions of weak amenability to fit this new framework, similarly to what was done in the article. After that, we would have to define the essential $C^*$-algebra of a \textbf{measured} groupoid opposed to a topological one and study the properties of this algebra relating it to the \textit{essential} weak amenability of the groupoid.

\end{document}